\documentclass[preprint,final,3p,times,12.8pt]{elsarticle}
\usepackage{mathrsfs}
\usepackage{amsthm,amsmath,amssymb}  
\usepackage{algorithm,algorithmicx,algpseudocode} 
\usepackage{array}
\usepackage{multirow}
\usepackage{bm}
\usepackage{hyperref}
\hypersetup{hypertex=true,
            colorlinks=true,
            linkcolor=blue,
            anchorcolor=blue,
            citecolor=blue}

\allowdisplaybreaks[4]
\newtheorem{theorem}{Theorem}[section]
\newtheorem{lemma}[theorem]{Lemma}
\newtheorem{assumption}[theorem]{Assumption}

\numberwithin{equation}{section}
\numberwithin{figure}{section}
\numberwithin{table}{section}



\journal{Applied Mathematics and Computation}

\begin{document}

\begin{frontmatter}

\title{An efficient solver based on low-rank approximation and Neumann matrix series for unsteady diffusion-type partial differential equations with random coefficients}

\author{Yujun Zhu\fnref{addr0}}
\ead{yujun\_zhu@hust.edu.cn}

\author{Min Li\fnref{addr0}}
\ead{minli@hust.edu.cn}

\author{Yulan Ning\fnref{addr0}}
\ead{ylning@hust.edu.cn}
		
\author{Ju Ming\fnref{addr0,addr1}\corref{cor1}} 
\ead{jming@hust.edu.cn}
\cortext[cor1]{Corresponding author}

\address[addr0]{School of Mathematics and Statistics, Huazhong University of Science and Technology, Wuhan, China}

\address[addr1]{Hubei Key Laboratory of Engineering Modeling and Scientific Computing, Huazhong University of Science and Technology, Wuhan, China}

\begin{abstract}
{
In this paper, we develop an efficient numerical solver for unsteady diffusion-type partial differential equations with random coefficients. A major computational challenge in such problems lies in repeatedly handling large-scale linear systems arising from spatial and temporal discretizations under uncertainty. To address this issue, we propose a novel generalized low-rank matrix approximation to represent the stochastic stiffness matrices, and approximate their inverses using the Neumann matrix series expansion. This approach transforms high-dimensional matrix inversion into a sequence of low-dimensional matrix multiplications. Therefore, the solver significantly reduces the computational cost and storage requirements while maintaining high numerical accuracy. The error analysis of the proposed solver is also provided. Finally, we apply the method to two classic uncertainty quantification problems: unsteady stochastic diffusion equations and the associated distributed optimal control problems. Numerical results demonstrate the feasibility and effectiveness of the proposed solver.
}
\end{abstract}


\end{frontmatter}

\section{Introduction}

Realistic simulations of complex systems governed by partial differential equations (PDEs) are fundamental in modern science and engineering. A prominent example is the diffusion equation, which models transport phenomena in heterogeneous media \cite{castrillon2021stochastic, gunzburger2011space}. It has been extensively used in a broad range of scientific and engineering applications, including chemical diffusion \cite{jha2016uncertainty, kuramoto2003chemical, mcrae1982numerical}, mass transport \cite{bird1956theory, cussler2009diffusion, flynn1974mass}, subsurface flow in porous media \cite{efendiev2003generalized, feng2019innovative, huyakorn2012computational,  miyan2015flow}, and acoustic energy propagation \cite{jing2007modified, su2019diffusion}. Nevertheless, in many practical scenarios, critical model parameters, e.g., permeability, are subject to considerable uncertainty,  arising from sparse data and incomplete physical characterization. This uncertainty limits the predictive capability of deterministic models and motivates the study of stochastic formulations.

To incorporate parameter uncertainty, deterministic diffusion equations are commonly reformulated as stochastic diffusion equations by treating the uncertain coefficients as random fields. In this framework, the solution itself becomes a random field defined over both physical and stochastic domains. The systematic study of such problems falls within the scope of uncertainty quantification (UQ) \cite{gunzburger2011error, peherstorfer2018survey, sullivan2015introduction}, which aims to characterize how uncertainties from the system inputs propagate through mathematical models and influence quantities of interest (QoIs). Therefore, developing robust and efficient UQ methodologies for stochastic diffusion equations is critical.

Realized by solving deterministic PDE for each independent and identically distributed (i.i.d.) realization of the random input, the Monte Carlo (MC) method \cite{fishman2013monte, metropolis1949monte} is a fundamental and widely used non-intrusive sampling technique in UQ. Its major advantage is that the convergence rate is independent of the dimensionality of the stochastic space. Despite its flexibility and ease of implementation, the MC method suffers from severe computational limitations when applied to time-dependent stochastic diffusion equations. This inefficiency arises from two factors. First, reliable statistical estimates of QoI require a large number of realizations, according to the Central Limit Theorem (CLT). Second, each realization needs to solve large-scale systems arising from spatial and temporal discretizations. When a fine spatial mesh is required, or when the stochastic partial differential equation (SPDE) is embedded in more complex settings such as optimal control, both computational cost and memory usage become prohibitively expensive. Consequently, the development of efficient solvers for these large-scale algebraic systems constitutes a key challenge in advancing UQ for time-dependent stochastic diffusion equations.

Dimensionality reduction offers a promising approach to alleviate this challenge, aiming to construct a more compact representation of the original high-dimensional data while maintaining the essential information. Many classical dimensionality reduction approaches rely on the vector space model \cite{turk1991eigenfaces, zhao2003face}, where samples are represented as vectors. Thus, the entire dataset is assembled into a data matrix and then projected onto a lower-dimensional subspace. Within this framework, low-rank matrix approximation via singular value decomposition (SVD) is one of the most well-known techniques, which provides optimal reconstruction under the Frobenius norm \cite{eckart1936approximation}. However, for high-dimensional matrices deriving from PDE discretizations, SVD imposes excessive demands in both computation and memory.

To overcome the limitations of the classical SVD method, the generalized low-rank approximation of matrices (GLRAM) was introduced in \cite{ye2004generalized}, which approximates a collection of matrices simultaneously by using a shared pair of two-sided orthogonal projection matrices. Numerous variants of GLRAM have been developed to enhance computational efficiency and robustness, which have been widely applied in computer vision and image analysis \cite{berry1995using, LI2022251, turk1991eigenfaces, zhao2003face}. For example, a non-iterative formulation was proposed to reduce computational cost in \cite{liu2006non}, and a simplified GLRAM structure was derived in \cite{lu2008simplified}. Robust extensions capable of handling sparse noise and outliers have also been investigated \cite{shi2015robust, zhao2016robust}. 

More recently, low-rank approximation techniques have been explored in the context of stochastic elliptic equations and the associated stochastic optimal control problems \cite{zhu2024low}. However, existing studies focus primarily on steady systems, and little attention has been paid to unsteady SPDEs. Moreover, inverting high-dimensional matrices remains a dominant computational bottleneck in time-dependent problems. These limitations motivate us to develop an unsteady SPDE solver that combines low-rank approximations with cost-effective matrix inversion schemes.

\begin{figure}[ht]
\centering
\includegraphics[width=1\textwidth]{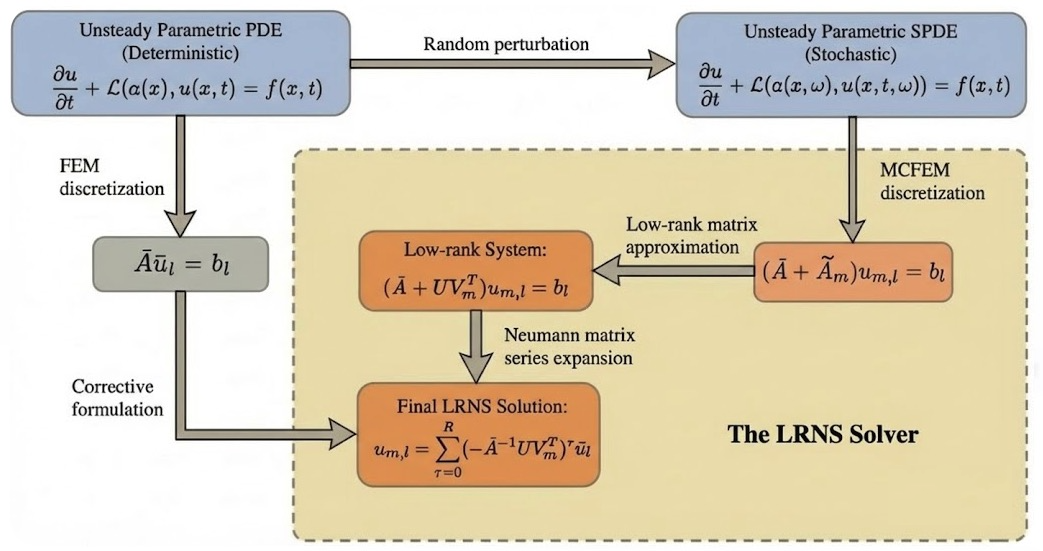}
\caption{The flowchart of the efficient solver.}
\label{fig1.1}
\end{figure}

In this paper, we develop an efficient solver for unsteady diffusion-type SPDEs. Specifically, we first adopt a novel generalized low-rank matrix approximation to represent the stochastic stiffness matrices, and then approximate their inverses using the Neumann matrix series expansion. Figure \ref{fig1.1} illustrates the schematic flowchart of our solver. The approach transforms high-dimensional matrix inversion into a sequence of low-dimensional matrix multiplications, thereby significantly reducing computational cost and storage requirements while maintaining high numerical accuracy. To demonstrate its effectiveness, the solver is employed to address two classical UQ problems: the unsteady diffusion equations with random permeability and the associated distributed optimal control problems. The numerical results confirm the feasibility and the efficiency of the proposed solver.

The remainder of this article is organized as follows. Section \ref{sec2} introduces the necessary function spaces and notations. Section \ref{sec3} presents a novel generalized low-rank matrix approximation and combines it with the Neumann series–based inversion strategy to construct an efficient solver for SPDEs. The error analysis is also provided in this section. In Section \ref{sec4}, we apply our solver to the unsteady stochastic diffusion equation and the corresponding stochastic optimal control problem. This section also reports the results of our numerical experiments. Conclusions and discussion can be found in Section \ref{sec5}.

\section{Preliminaries}\label{sec2}

We begin by recalling some necessary function spaces and notations. Throughout this article, we use the standard notations for Sobolev spaces \citep{adams2003sobolev}. Let $L^p(\mathcal{D}),  1 \leq p \leq \infty$, denote usual Lebesgue space on smooth domain $\mathcal{D} \in \mathbb{R}^n$; $\Vert \cdot \Vert = \Vert \cdot \Vert_{L^2(\mathcal{D})}$ denote the $L^2(\mathcal{D})$-norm induced by the inner product $\langle f,g \rangle = \displaystyle \int_\mathcal{D} fg \,\mathrm{d}\mathcal{D}, \, \forall f,g\in L^2(\mathcal{D})$. $H^r(\mathcal{D})$ is a Sobolev space for all real numbers $r$ with norms $\Vert y \Vert_{H^r(\mathcal{D})}$ and semi-norms $\vert y \vert_{H^r(\mathcal{D})}$, where

\begin{equation} \label{eq2.1}
\Vert y \Vert^2_{H^r(\mathcal{D})} = \sum_{\vert \bm{\alpha} \vert \leq r} \Vert \frac{\partial^{\vert\bm{\alpha}\vert} y}{\partial x^{\bm{\alpha}}} \Vert^2_{L^2(\mathcal{D})},
\end{equation}

\noindent and

\begin{equation} \label{eq2.2}
\vert y \vert^2_{H^r(\mathcal{D})} = \sum_{\vert \bm{\alpha} \vert = r} \Vert \frac{\partial^{\vert\bm{\alpha}\vert} y}{\partial x^{\bm{\alpha}}} \Vert^2_{L^2(\mathcal{D})}.
\end{equation}

\noindent Here $\alpha$ is a multi-index with non-negative integer components $\{\alpha_i\}$, and $\bm{\alpha} = \sum_i \alpha_i$.

Sobolev spaces have forms

\begin{equation} \label{eq2.3}
H^1(\mathcal{D}) = \Big\{ \,y \in L^2(\mathcal{D}), \partial_{x_i} y \in L^2(\mathcal{D}), i = 1,...,n\,\Big\},
\end{equation}

\noindent and

\begin{equation} \label{eq2.4}
H_0^1(\mathcal{D}) = \Big\{ \, y \in H^1(\mathcal{D}), y \vert_{\partial \mathcal{D}} = 0\,\Big\}.
\end{equation}

\noindent Clearly, $H_0^1(\mathcal{D})$ is a subspace of $H^1(\mathcal{D})$ and it has the dual space $H^{-1}(\mathcal{D})$. Then we define the Hilbert space $L^2(0,T; H^r(\mathcal{D})) = \left\{ y(t) \in H^r(\mathcal{D}), \int_\mathcal{D} \Vert y \Vert_{H^r(\mathcal{D})} \mathrm{d}t  <  \infty \right\}$ with the norm $\Vert y \Vert_{L^2(0,T; H^r(\mathcal{D}))} = \left(\int_0^T \Vert y \Vert_{H^r(\mathcal{D})} \mathrm{d}t\right)^{\frac{1}{2}}$.

Let $(\Omega, \mathscr{F}, \mathbb{P})$ be a complete probability space, where $\Omega$ denotes a set of random events, $\mathscr{F}$ is the minimal $\sigma$-algebra of outcomes, and $\mathbb{P}: \mathscr{F} \rightarrow [0,1]$ is a complete probability measure, respectively. If $X$ is a real random variable in $(\Omega, \mathscr{F}, \mathbb{P})$, then its expectation is given by

\begin{equation} \label{eq2.5}
\mathbb{E}[X] = \int_\Omega X(\omega) \mathbb{P}(d \omega) = \int_{\mathbb{R}^n} x \rho(x) dx,
\end{equation}

\noindent where $\rho$ is the joint probability distribution measure for $X$, and it is defined on a Borel set $\mathcal{B}$ of $\mathbb{R}$, i.e., $\rho(\mathcal{B}) = \mathbb{P}(X^{-1}(\mathcal{B}))$.

Finally, we define the stochastic Sobolev space as

\begin{equation} \label{eq2.6}
L^2(\Omega;L^2(0,T; H^r(\mathcal{D}))) = \left\{ y: \Omega \rightarrow L^2(0,T; H^r(\mathcal{D})), \mathbb{E}\left[\Vert y \Vert^2_{L^2(0,T; H^r(\mathcal{D}))}\right]< \infty\right\},
\end{equation}

\noindent equipped with the norm

\begin{equation}\label{eq2.7}
\Vert y \Vert^2_{L^2(\Omega;L^2(0,T; H^r(\mathcal{D})))} = \mathbb{E}\left[ \Vert y \Vert_{L^2(0,T; H^r(\mathcal{D}))}\right].
\end{equation}

The stochastic Sobolev space $L^2(\Omega;L^2(0,T; H^r(\mathcal{D})))$ is a Hilbert space, and it is isomorphic to the tensor product space $L^2(\Omega) \otimes L^2(0,T; H^r(\mathcal{D}))$. We define the space-time cylinder $\mathcal{Q} = (0, T] \times \mathcal{D}$ for a given $T > 0$. For simplicity, we set 
\begin{equation} \label{eq2.8}
\mathcal{H}^r(\mathcal{Q}) = L^2(\Omega;L^2(0,T; H^r(\mathcal{D}))), \enspace\mathcal{H}_0^r(D) = L^2(\Omega;L^2(0,T; H_0^r(\mathcal{D}))), \enspace\mathcal{L}^2(D) = L^2(\Omega;L^2(0,T; L^2(\mathcal{D}))).
\end{equation}

\section{An efficient solver for perturbed SPDEs via low-rank approximation and Neumann matrix series (LRNS) } \label{sec3}

The primary challenge in solving SPDEs lies in the need to handle a large number of algebraic systems, which often results in prohibitive computational and storage costs. To alleviate this issue, we introduce a novel low-rank matrix approximation approach based on randomized singular value decomposition (RSVD) in this section, which can significantly reduce the memory requirement. By integrating this low-rank representation with the Neumann matrix series, we develop an efficient solver, i.e., LRNS solver, for perturbed SPDEs and carry out its error analysis.

\subsection{Problem statement and perturbation splitting} \label{subsec3.1}

Consider a PDE operator $\mathcal{L}$ and a random field $a(\bm{x},\omega): \mathcal{D} \times \Omega \rightarrow \mathbb{R}$. For a given $T > 0$ and a deterministic force term $f(\bm{x},t) \in \mathcal{Q}:=\mathcal{D} \times (0,T]$, we study the following unsteady parametric SPDE: find $u(\bm{x}, t,\omega): \mathcal{Q} \times \Omega \rightarrow \mathbb{R}$ such that

\begin{equation} \label{eq3.1}
\frac{\partial u}{\partial t} + \mathcal{L}(a(\bm{x},\omega), u(\bm{x}, t,\omega)) = f(\bm{x},t), \enspace a.e. \enspace in \enspace \mathcal{Q}, 
\end{equation}

\noindent on a bounded, Lipschitz domain $\mathcal{D} \subset \mathbb{R}^d$, $d = 1,2,3$ and a complete probability space $(\Omega, \mathscr{F}, \mathbb{P})$, subject to suitable boundary conditions. The finite-dimensional random variable $\omega$ in $a(\mathbf{x}, \omega)$ admits a joint probability density function (PDF) $\rho: \Omega \rightarrow \mathbb{R}^+$, with $\rho \in L^\infty(\Omega)$. 

The random parameter $a(\bm{x},\omega)$ in the time-dependent stochastic systems \eqref{eq3.1} is generally assumed to have uniform positive upper and lower bounds. A commonly used representation in heterogeneous porous media is the log-normal formulation \cite{babuvska2007stochastic, ganis2008stochastic, dostert2008multiscale}, i.e., $a(\bm{x},\omega) = \exp\left(\eta(\bm{x},\omega)\right)$, where the field $\eta(\bm{x},\omega)$ follows a Gaussian distribution. Although the formulation guarantees positivity, the resulting log-normal distribution does not always align with statistical features observed in field data. In many practical settings, such as regression analysis and geo-statistical inversion, perturbations in model parameters are more naturally captured by an additive Gaussian structure \cite{ali2017multilevel, barth2011multi, betz2014numerical, ghanem1998stochastic}. To remain consistent with statistical practices, we adopt the following additive perturbation formulation for the random parameter in \eqref{eq3.1}:

\begin{equation} \label{eq3.2}
a(\bm{x},\omega) = \bar{a}(\bm{x}) + \widetilde{a}(\bm{x},\omega),
\end{equation}

\noindent where $\widetilde{a}(\bm{x},\omega)$ represents a random perturbation around the deterministic mean field $\bar{a}(\bm{x})$. The mean $\bar{a}(\bm{x})$ is the dominant component of the permeability field and ensures its positivity and boundedness. The perturbation setting can accurately depict many real-world scenarios in which a well-understood physical process is subject to random variations. It has been widely adopted in industrial and engineering applications, including surface and subsurface flows \cite{freeze1972role, furman2008modeling, huyakorn2012computational, sudicky2008simulating}, heat transfer in composite materials \cite{kalamkarov2014three, sakata2015successive, shao2017review}, and structural mechanics \cite{chen2024uncertainty, kaminski2013stochastic, lowdin1951note, nayfeh2024perturbation}. For instance, in composite materials such as fiberglass or reinforced concrete, the thermal conductivity is determined by bulk properties and volume fractions of the components. However, it inevitably exhibits slight random deviations due to manufacturing imperfections, local density variations, and microscopic defects. This situation naturally falls within the scope of stochastic perturbation theory.

To quantify uncertainties in the physical system \eqref{eq3.1}, our goal is to obtain the QoI that provides a probabilistic description of the system's behavior, rather than a single realization of $u$. Specifically, we aim to compute the multi-dimensional statistical expectation \citep{smith2013uncertainty}:

\begin{equation} \label{eq3.3}
\mathbb{E}[u](\bm{x},t) = \int_\Omega u(\bm{x}, t,\omega) \rho(\omega) d\omega, \quad in \enspace \mathcal{Q} \times \Omega. 
\end{equation}

It is generally infeasible to obtain analytic solutions of the SPDE in \eqref{eq3.1}. A natural choice for numerical implementation is the Monte Carlo finite element method (MC-FEM). It is a robust, non-intrusive approach that can effectively mitigate the curse of dimensionality. We generate $M$ i.i.d. realizations $\widetilde{a}_m(x) := \widetilde{a}(x,\omega_m)$ and discretize $M$ independent deterministic unsteady PDEs by the finite element (FE) method in space and the finite difference (FD) scheme in time. Let $N$ and $L$ denote the number of spatial degrees of freedom and sub-intervals of duration respectively. Then, the resulting algebraic systems for the solution vector $\bm{u}_{m,l} \in \mathbb{R}^N$ is given by

\begin{equation} \label{eq3.4}
\mathbf{A}_m \bm{u}_{m,l} = \bm{b}_l,\quad m=1,2,\cdots,M, \:l = 1,2,\cdots, L, 
\end{equation}

\noindent where each stiffness matrix $\mathbf{A}_m$ can be decomposed into the deterministic stiffness matrix $\bar{\mathbf{A}}$ derived from $\bar{a}(\bm{x})$, and the sample-dependent perturbation matrix $\widetilde{\mathbf{A}}_m$ arising from $\widetilde{a}(\bm{x},\omega)$:
 
\begin{equation}\label{eq3.5}
\mathbf{A}_m=\bar{\mathbf{A}} \: + \: \widetilde{\mathbf{A}}_m,\quad m=1,2,\cdots,M.
\end{equation}

From the perspective of signal processing and deep learning \citep{hemanth2017deep}, the pre-processing of subtracting the statistical mean is also known as data normalization, which removes the common parts in the MC samples, highlights the individual differences, thereby improveing the generalization ability of our algorithm. 

The QoI in \eqref{eq3.3} is then approximated as

\begin{equation} \label{eq3.6}
\mathbb{E}[u](\bm{x},t)  \: \approx \:  \left(\bm{u}_{0},\cdots,\bm{u}_{L} \right) \: = \: \left(\bm{u}_{0},\frac{1}{M} \sum_{m=1}^M \bm{u}_{m,1},\cdots,\frac{1}{M} \sum_{m=1}^M \bm{u}_{m,L}\right),  
\end{equation}

\noindent where the initial state $\bm{u}_{0}$ is determined by the initial conditions from the SPDE system.

The standard MC-FEM requires solving large linear systems involving a collection of matrices $\{\widetilde{\mathbf{A}}_m\}_{m=1}^M$ at each time step. In many UQ problems, both the number of MC samples $M$ and the number of degrees of freedom $N$ must be large to achieve the desired precision, which inevitably leads to massive computational cost and memory requirements for storing and inverting $M \times L$ distinct sparse but large stiffness matrices. These challenges motivate us to develop a high-performance SPDE solver that can effectively mitigate the unaffordable computational and space complexities. Therefore, we introduce a novel low-rank approximation for $\widetilde{\mathbf{A}}_m$ in Section \ref{subsec3.2}, and present an efficient matrix series inversion technique in Section \ref{subsec3.3}.

\subsection{RSVD-based low-rank approximation of stiffness matrices} \label{subsec3.2}

The low-rank approximation (LRA) is a fundamental and popular dimensionality reduction technique. Its objective is to approximate high-dimensional matrices with lower-rank representations that retain the essential intrinsic information. By operating on these more compact surrogate matrices, LRA effectively alleviates the storage and computational burdens while maintaining approximately identical numerical accuracy to that achieved with the full matrices.

Mathematically, given a matrix $\mathbf{B} \in \mathbb{R}^{N \times N}$, its optimal rank-$k$ approximation ($k \le N$) under the Frobenius norm is defined by the following minimization problem:

\begin{equation} \label{eq3.7}
	 \mathop{\min}\limits_{\text{rank}(\mathbf{B}^*) \; = \; k} \enspace \Vert \mathbf{B} - \mathbf{B}^* \Vert_F,
\end{equation}

\noindent where $ \Vert \mathbf{B} - \mathbf{B}^* \Vert_F$ represents the reconstruction error between the original and approximate matrices.

According to the famous Eckart-Young-Mirsky theorem, the rank-constrained minimization problem in \eqref{eq3.7} admits a closed-form solution through the truncated SVD.

\begin{theorem}[\cite{eckart1936approximation}]\label{th3.1} 

Let $\mathbf{B} = \mathbf{U} \Sigma \mathbf{V}^T$ be the SVD of $\mathbf{B}$, where $\mathbf{U},\mathbf{V}$ are orthogonal matrices and $\Sigma$ is a diagonal matrix with the singular values of $\mathbf{B}$ in descending order. Then, the optimal rank-$k$ approximation of $\mathbf{B}$ under the Frobenius norm is given by
$$\mathbf{B}^* \:= \:\mathbf{U}_k \Sigma_k \mathbf{V}_k^T,$$ 
\noindent where $\mathbf{U}_k, \mathbf{V}_k \in \mathbb{R}^{N \times k}$  and $\Sigma_k \in \mathbb{R}^{k \times k}$. The minimizer $\mathbf{B}^*$ is unique if and only if $\sigma_{k+1} \ne \sigma_k$.

\end{theorem}

Although SVD provides the optimal low-rank approximation for a single matrix, it becomes impractical to deal with a large collection of matrices, as applying SVD independently to each matrix often results in prohibitive computational and memory costs. To overcome the limitation, the generalized low-rank approximation of matrices framework \cite{ye2004generalized} was introduced. For a sequence of matrices $\{\mathbf{B}_m\}_{m=1}^M$, GLRAM seeks a shared pair of orthonormal matrices $\mathbf{U},\mathbf{V} \in \mathbb{R}^{N \times k}$ and a set of reduced representations $\{\mathbf{S}_m\}_{m=1}^M \in \mathbb{R}^{k \times k}$ by solving the following optimization formulation:

\begin{equation} \label{eq3.8}
\mathop{\min}\limits_{\mathbf{U}^T \mathbf{U} = \mathbf{V}^T \mathbf{V} = \mathbf{I}_k} \quad \sum_{m=1}^M \Vert \mathbf{B}_m - \mathbf{U} \mathbf{S}_m \mathbf{V}^T \Vert_F^2.
\end{equation}

Due to its significant enhancement in computational efficiency and data compression performance, GLRAM extends the applicability of SVD-based data compression techniques to the multi-matrix settings. These appealing properties motivate us to incorporate the concept of "generalized" to solve the complex algebraic systems in \eqref{eq3.4}. However, the classical GLRAM method still faces several challenges. First, it relies on an iterative structure that requires two eigen-decompositions of large matrices at each iteration \cite{ye2004generalized}, which makes the algorithm still computationally demanding. Second, GLRAM has been primarily used in computer vision and signal processing, making it nontrivial to achieve the reconstruction precision required for high-fidelity PDE simulations. 

To address these issues, we propose a novel GLRAM strategy that achieves a better trade-off among dimensionality reduction effectiveness, computational efficiency, and reconstruction accuracy. Unlike the conventional SVD and GLRAM approaches, it employs a different low-rank matrix representation that significantly reduces both computational and space complexities while maintaining high fidelity in matrix reconstruction. Moreover, the novel data compression method utilizes the RSVD technique to further accelerate the matrix computations. As a result, our approach exhibits robust performance in handling the perturbed linear systems in \eqref{eq3.4}. 

We begin by introducing a prescribed data compression ratio $\tau \in [0,1]$ that determines the reduced dimension as below,

\begin{equation}
k \: = \: \lceil \tau N \rceil \: \in \: [0,N].
\nonumber
\end{equation}

\noindent Our goal is to find the optimal solutions $\mathbf{U}, \{\mathbf{V}_m\}_{m=1}^M \in \mathbb{R}^{N \times k}$ by solving the following optimization problem:

\begin{equation} \label{eq3.9}
	\mathop{\min}\limits_{ \substack{\mathbf{U},\mathbf{V}_m \in \mathbb{R}^{N \times k} \\ \mathbf{U}^T \mathbf{U} = \mathbf{I}_k}} \quad \sum_{m=1}^M \Vert \mathbf{B}_m - \mathbf{U} \mathbf{V}_m^T \Vert_F^2.
\end{equation}

Then, the following theorems provide a closed-form, non-iterative formulation to obtain the optimal solutions $\mathbf{U}, \{\mathbf{V}_m\}_{m=1}^M$ in \eqref{eq3.9}. 

\begin{theorem}[\cite{zhu2024low}] \label{th3.2}
	
	Let $\mathbf{U}$ and $\{\mathbf{V}_m\}_{m=1}^M$ be the optimal solution of the minimization problem in \eqref{eq3.9}. Then for each $m = 1,2,\cdots,M$, 

    \begin{equation} \label{eq3.10}
	   \mathbf{V}_m \:= \:\mathbf{B}_m^T \mathbf{U}.
    \end{equation}
	
\end{theorem}

\begin{theorem}[\cite{zhu2024low}]\label{th3.3}
	
	Let $\mathbf{U}$ be the optimal solution to the minimization problem in \eqref{eq3.9}. Then, $\mathbf{U}$ consists of the $k$ eigenvectors corresponding to the largest $k$ eigenvalues of
	
	\begin{equation} \label{eq3.11}
		\mathbf{N} =  \sum_{m=1}^M \mathbf{B}_m \mathbf{B}_m^T.
	\end{equation}

\end{theorem}

Although we derive a straightforward non-iterative algorithm, Theorem \ref{th3.3} indicates that it requires computing the eigen-decomposition of the dense matrix $\mathbf{N} \in \mathbb{R}^{N \times N}$, which becomes prohibitively expensive when the degree of freedom $N$ is large. To alleviate this burden, we integrate the RSVD technique \cite{halko2011finding} into the analytical framework. 

RSVD is a well-established randomized technique that provides inexpensive yet accurate low-rank approximations. Its core idea is to use randomization to identify a low-dimensional subspace that captures the dominant spectral features of a large matrix. By projecting the matrix onto the randomly generated lower-dimensional subspace, RSVD avoids expensive calculations on the full matrix, thereby significantly reducing the overall computational cost while maintaining high approximation precision. For the matrix $\mathbf{N} \in \mathbb{R}^{N \times N}$, the RSVD implementation involves two main stages:

\begin{enumerate}
    \item \textit{Subspace Identification}: Initialize a random Gaussian matrix $\mathbf{P} \in \mathbb{R}^{N \times k}$ and form the sketch $\mathbf{Z} = \mathbf{N}\mathbf{P}$, whose columns are linear combinations of the columns of $\mathbf{N}$. Perform a QR decomposition $\mathbf{Z} = \mathbf{Q}\mathbf{R}$, where $\mathbf{Q} \in \mathbb{R}^{N \times k}$ provides an orthonormal basis for the range of $\mathbf{Z}$. The columns of $\mathbf{Q}$ probabilistically span a $k$-dimensional subspace that well approximates the dominant eigen-subspace of $\mathbf{N}$.

    \item \textit{Projection}: Project the original matrix $\mathbf{N}$ onto the low-dimensional subspace spanned by $\mathbf{Q}$ to obtain $\mathbf{Y} = \mathbf{Q}^T \mathbf{N}$. Compute its SVD $\mathbf{Y} = \mathbf{U}_Y \Sigma_Y \mathbf{V}_Y^T$, where $\mathbf{U}_Y,\Sigma_Y \in \mathbb{R}^{k \times k}$ and $\mathbf{V}_Y \in \mathbb{R}^{N \times k}$. The approximate eigen-vectors of $\mathbf{N}$ are recovered by $\mathbf{U} = \mathbf{Q} \mathbf{U}_Y$, and the rank-$k$ RSVD formulation of $\mathbf{N}$ is then given by $\mathbf{N} = \mathbf{U} \Sigma_Y \mathbf{V}_Y^T$.
\end{enumerate}

The classical SVD of $\mathbf{N} \in \mathbb{R}^{N \times N}$ requires $\mathcal{O}(N^3)$ operations, while RSVD only requires $\mathcal{O}(N^2k)$ operations \cite{halko2011finding}. The following theorem also ensures the accuracy of RSVD. 

\begin{theorem}[\cite{halko2011finding}]\label{th3.4} 

Let the eigenvalues of $\mathbf{N}$ are $\lambda_1(\mathbf{N}) \ge \lambda_2(\mathbf{N})\ge \cdots \ge \lambda_N(\mathbf{N})$, then expected approximation error satisfies
$$\mathbb{E}\Vert\mathbf{N} - \mathbf{Q}\mathbf{Q}^T\mathbf{N}\Vert_F \:\le\: \sqrt{1+k} \,\left(\sum_{i=k+1}^N \lambda_i(\mathbf{N})\right)^{\frac{1}{2}}.$$
\end{theorem}

Finally, we can summarize the complete procedure for determining the optimal solution $\mathbf{U}$ and $\{\mathbf{V}_m\}_{m=1}^M$ of \eqref{eq3.9} in Algorithm \ref{alg1}.

\begin{algorithm}[htbp]
\caption{A novel RSVD-based generalized low-rank approximation method of matrices} \label{alg1}
\begin{algorithmic}[1]
\Require
Matrices $\{\mathbf{B}_m\}_{m=1}^M$, and the data compression ratio $\tau$
\Ensure
Matrices $\mathbf{U}$ and $\{\mathbf{V}_m\}_{m=1}^M$
		
\State Compute the reduced rank $k = \lceil \tau N \rceil$. 
		
\State Form the matrix $\mathbf{N}$ in \eqref{eq3.11}.

\State Generate a random $N \times k$ sketch matrix and determine the eigenvectors of $\mathbf{N}$ via RSVD.
		
\State Construct the matrix $\mathbf{U}$ by taking $k$ eigenvectors of $\mathbf{N}$ that correspond to its largest $k$ eigenvalues.

\State Compute the matrices $\mathbf{V}_m$ by \eqref{eq3.10} for each $m = 1,\cdots,M$.
		
\State \Return the optimal solution $\mathbf{U}$ and $\{\mathbf{V}\}_{m=1}^M$ of \eqref{eq3.9}.
		
\end{algorithmic}
\end{algorithm}

\begin{figure}[!htbp]
\centering
\includegraphics[width=1\textwidth]{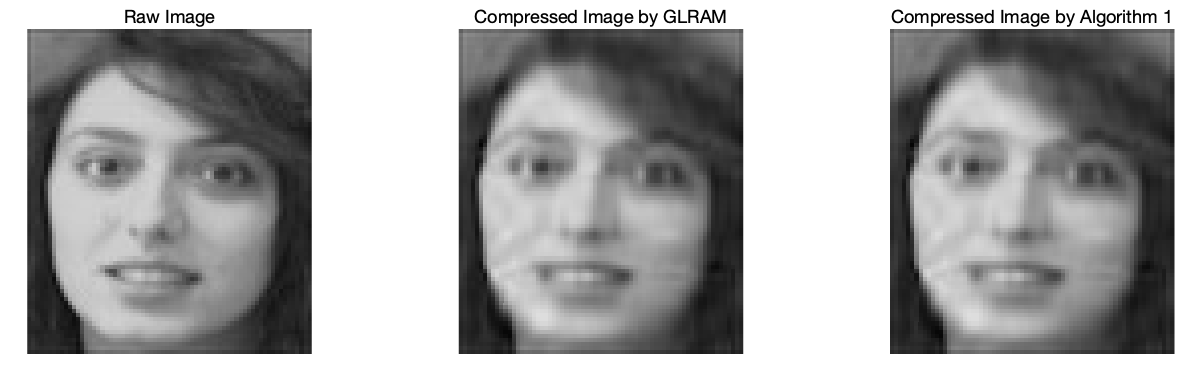}
\caption{Raw image from LFW dataset (left), the images compressed by conventional GLRAM (middle) and Algorithm \ref{alg1} (right) using the data compression ratio $\tau=10\%$.}
\label{fig3.1}
\end{figure}

To highlight the advantages of the proposed approach, Table \ref{tab3.1} compares the computational and storage complexities among the standard SVD, conventional GLRAM, and Algorithm \ref{alg1}. The results clearly demonstrate that our method achieves notable improvements in both runtime and memory usage, particularly when dealing with large collections of matrices. The last row of Table \ref{tab3.1} indicates that the improvement scales roughly linearly with the data compression ratio $\tau$. Selecting a relatively small $\tau$ accelerates the matrix computations and reduces storage demand. Nonetheless, pushing $\tau$ too far may lose critical information from the original matrices. The trade-off between computational efficiency and approximation accuracy will be analyzed in the subsequent sections. Figure \ref{fig3.1} also demonstrates the nice compression performance of Algorithm \ref{alg1}, which could well capture the main features of the face image while requiring small storage. Here, the image is randomly chosen from the LFW dataset \cite{huang:inria-00321923}.

\begin{table}[!ht] 
\centering 
\caption{Comparison among the traditional SVD, standard GLRAM, and Algorithm \ref{alg1}. Here $M$, $N$, $k$, and $I$ denote the sample size, original dimension, reduced dimension, and iteration count in standard GLRAM, respectively.}
\begin{tabular}{ccc} \hline 
Method & Space & Time  \\ \hline
Standard SVD & $MN^2$  & $\mathcal{O}(MN^3)$\\
Standard GLRAM & $2Nk + Mk^2$ & $\mathcal{O}(2IMN^2 k+ 2IN^3)$ \\
Algorithm \ref{alg1} & $Nk + MNk$ & $\mathcal{O}(M N^2 k + Nk^2 + k^3)$\\ 
Ratio between SVD and Algorithm \ref{alg1} & $\mathcal{O}(\tau)$ & $\mathcal{O}(\tau)$ \\ \hline
\end{tabular}
\label{tab3.1}
\end{table}

\subsection{Approximation of inverse matrices via Neumann series} \label{subsec3.3}

The algebraic systems in \eqref{eq3.4}, derived from the standard MC-FEM and FD schemes, are computationally intensive and memory-consuming. To overcome the challenge, we combine the novel RSVD-based LRA technique introduced in Section \ref{subsec3.2} with the Neumann matrix series expansion, which allows us to efficiently approximate the required matrix inverses while maintaining accuracy.

Recall that each full stiffness matrix $\mathbf{A}_m$ is decomposed into a constant deterministic matrix $\bar{\mathbf{A}}$ and a sample-dependent perturbation matrix $\widetilde{\mathbf{A}}_m$ that captures the randomness in the system. Under the perturbation assumption in \eqref{eq3.2} and spatially correlated nature of $\widetilde{a}$, Algorithm \ref{alg1} provides the following low-rank approximation of $\widetilde{\mathbf{A}}_m$:

$$\widetilde{\mathbf{A}}_m \: \approx \: \widetilde{\mathbf{A}}_m^* \: := \: \mathbf{U} \mathbf{V}_m^T, \quad m = 1,\cdots,M,$$ 

\noindent where $\mathbf{U}, \mathbf{V}_m \in \mathbb{R}^{N \times k}$. Then the linear formulation in \eqref{eq3.4} is transformed to:

\begin{equation} \label{eq3.12}
\begin{aligned}
\bm{u}_{m,l} \: &=\: \left(\bar{\mathbf{A}} + \widetilde{\mathbf{A}}_m^*\right)^{-1} \bm{b}_l\\
&= \:\left[\bar{\mathbf{A}} \left( \mathbf{I}_N + \bar{\mathbf{A}}^{-1}\widetilde{\mathbf{A}}_m^*\right)\right]^{-1} \bm{b}_l\\
&= \:\left( \mathbf{I}_N + \bar{\mathbf{A}}^{-1}\widetilde{\mathbf{A}}_m^*\right)^{-1} \bar{\mathbf{A}}^{-1}\bm{b}_l\\
&= \:\left( \mathbf{I}_N + \bar{\mathbf{A}}^{-1}\widetilde{\mathbf{A}}_m^*\right)^{-1} \bar{\bm{u}}_l, 
\end{aligned}
\end{equation}

\noindent for $m=1,2,\cdots,M, \;l = 1,2,\cdots, L$. Here $\mathbf{I}_N$ represents the $N$-dimensional identity matrix and $\bar{\bm{u}}_l = \bar{\mathbf{A}}^{-1}\bm{b}_l$ is the unperturbed numerical solution at the $l$-th time step.

The main challenge is now switched to approximating $\left( \mathbf{I}_N + \bar{\mathbf{A}}^{-1}\widetilde{\mathbf{A}}_m^* \right)^{-1}$, which can also be interpreted as a low-rank correction $\bar{\mathbf{A}}^{-1}\widetilde{\mathbf{A}}_m^*$ to the identity matrix $\mathbf{I}_N$. The special corrective structure inspires us to employ the Neumann series, also known as the geometric series for matrices, which follows the following lemma.

\begin{lemma} \label{le3.5}
Let $\mathbf{B} \in \mathbb{R}^{N\times N}$ and suppose that its spectral radius satisfies $\rho(\mathbf{B}) < 1$. Then $\mathbf{I}_N - \mathbf{B}$ is non-singular and the Neumann series $\displaystyle \sum_{r=0}^\infty \mathbf{B}^r$ converges to the inverse $\left(\mathbf{I}_N - \mathbf{B}\right)^{-1}$.
\end{lemma}

According to Lemma \ref{le3.5}, the inverse of $\mathbf{I}_N + \bar{\mathbf{A}}^{-1}\widetilde{\mathbf{A}}_m^*$ is approximated as follows

\begin{equation} \label{eq3.13}
\left(\mathbf{I}_N + \bar{\mathbf{A}}^{-1}\widetilde{\mathbf{A}}_m^* \right)^{-1} \:=\: \sum_{r=0}^\infty \left(-\bar{\mathbf{A}}^{-1}\widetilde{\mathbf{A}}_m^*\right)^{r},
\end{equation}

\noindent whose convergence condition is that the spectral radius $\rho(\bar{\mathbf{A}}^{-1}\widetilde{\mathbf{A}}_m^*) < 1$. In practice, since $\widetilde{a}(\bm{x},\omega)$ represents only a small fluctuation around the mean coefficient $\bar{a}(\bm{x})$, the norm $\Vert\widetilde{\mathbf{A}}_m^*\Vert$ remains sufficiently small. Consequently, the spectral radius can be bounded as

$$\rho(\bar{\mathbf{A}}^{-1}\widetilde{\mathbf{A}}_m^*) \:=\: \Vert\bar{\mathbf{A}}^{-1}\widetilde{\mathbf{A}}_m^*\Vert  \:\le\: \Vert\bar{\mathbf{A}}^{-1}\Vert\Vert\widetilde{\mathbf{A}}_m^*\Vert \:<\: 1,$$

\noindent which guarantees the convergence of \eqref{eq3.13}. Hence, employing the Neumann series expansion to approximate the matrix inverse is well-suited for time-dependent PDEs involving perturbations in the coefficient matrices.
 
For computational efficiency, the infinite Neumann series is truncated after a finite number of terms $R$:

$$\left(\mathbf{I}_N + \bar{\mathbf{A}}^{-1}\widetilde{\mathbf{A}}_m^* \right)^{-1} \:\approx\: \sum_{r=0}^R \left(-\bar{\mathbf{A}}^{-1}\widetilde{\mathbf{A}}_m^*\right)^{r}.$$

\noindent Note that the terms $\left(-\bar{\mathbf{A}}^{-1}\widetilde{\mathbf{A}}_m^*\right)^{r}$ in the series can be computed in parallel, leading to higher throughput. According to \eqref{eq3.13}, the truncation error is proportional to $\Vert(-\bar{\mathbf{A}}^{-1}\widetilde{\mathbf{A}}_m^*)^{R+1}\Vert$. As a result, the accuracy of this approximation grows exponentially as the truncation index $R$ increases. 

Finally, by substituting this approximation into \eqref{eq3.12}, the numerical solutions of the algebraic systems in \eqref{eq3.4} are given by

\begin{equation}\label{eq3.14}
\mathbf{u}_{m,l} = \sum_{r=0}^R \left(-\bar{\mathbf{A}}^{-1}\mathbf{U}\mathbf{V}_m^T\right)^{r} \bar{\bm{u}}_l,\quad m=1,2,\cdots,M, \:l = 1,2,\cdots, L,
\end{equation}

\noindent which also indicates that $\mathbf{u}_{m,l}$ can be interpreted as a series-based correction to the deterministic unperturbed state $\bar{\bm{u}}_l = \bar{\mathbf{A}}^{-1} \bm{b}_l$. 

Therefore, we develop an efficient solver based on the low-rank matrix approximation and Neumann matrix series expansion, referred to as the LRNS solver, for the linear systems in \eqref{eq3.4} arising from the full discretization of parameter-dependent PDEs with perturbations. The solver integrates the novel RSVD-based LRA framework in Algorithm \ref{alg1} with the series-based correction strategy in \eqref{eq3.14}, which effectively transforms expensive matrix inversions into a series of low-cost corrective operations applied to the precomputed deterministic solution. The entire procedure avoids the conventional direct inversion of each full matrix $\mathbf{A}_m$ and relies solely on matrix–vector multiplications and low-rank updates, thereby significantly reducing both computational and memory complexities while preserving satisfactory numerical accuracy. The pseudo-code for the LRNS procedure is demonstrated in Algorithm \ref{alg2}.

\begin{algorithm}[h]
\caption{The LRNS solver for the algebraic systems originated from the full discretization of parameter-dependent PDEs with perturbations} 
\label{alg2}
\begin{algorithmic}[1]

\Require
Data compression ratio $\tau$, MC sample size $M$, numbers of spatial and temporal discretization for SPDEs $N$ and $L$, Neumann series truncation index $R$, deterministic stiffness matrix $\bar{\mathbf{A}}$, perturbed stiffness matrices $\{\widetilde{\mathbf{A}}_m\}_{m=1}^M$, and load vectors for each time step $\{\bm{b}_l\}_{l=0}^L$.

\Ensure
Approximation of QoI $\mu$.

\State Compute the deterministic solution $\bar{\bm{u}}_l = \bar{\mathbf{A}}^{-1} \bm{b}_l$ for $l =1,2,\cdots,L$ and set the initial state $\bm{u}^0$.

\State Determine the reduced dimension $k = \lceil \tau N \rceil$ and obtain the optimal rank-$k$ approximation $\mathbf{U}, \{\mathbf{V}_m\}_{m=1}^M$ of $\{\widetilde{\mathbf{A}}_m\}_{m=1}^M$ using Algorithm \ref{alg1}.

\For{$m = 1, \ldots, M$}  
\For{$l = 1, \cdots, L$}
    \State Compute the sample solution $\bm{u}_{m,l}$ by the corrective formulation \eqref{eq3.14}.
\EndFor
\EndFor

\State \Return the estimated QoI $\mu$ by \eqref{eq3.6}.

\end{algorithmic}
\end{algorithm}

\subsection{Error analysis}\label{subsec3.4}

This subsection provides the error analysis of the LRNS solver in Algorithm \ref{alg2} by separately studying the errors
from low-rank matrix representation and inverse approximation via Neumann series, respectively. We begin by introducing a common metric for assessing the LRA performance, the Root Mean Square Reconstruction Error (RMSRE):

\begin{equation} \label{eq3.15}
\text{RMSRE} \: := \: \sqrt{\frac{1}{M} \sum_{m=1}^M \Vert \widetilde{\mathbf{A}}_m - \mathbf{U} \mathbf{V}_m^T \Vert_F^2},
\end{equation}

\noindent which quantifies the overall discrepancy between the original full matrices and their low-rank approximation. A smaller RMSRE generally indicates a more accurate approximation to the collection $\{\widetilde{\mathbf{A}}_m\}_{m=1}^M$.

To examine the RMSRE value achieved by our novel low-rank approximation method in Algorithm \ref{alg1}, we first look at the following theorem.

\begin{theorem} [\cite{zhu2024low}] \label{th3.7}
Let $\mathbf{U}, \{\mathbf{V}_m\}_{m=1}^M \in \mathbb{R}^{N \times k}$ be the optimal solutions of the minimization problem in \eqref{eq3.9}. Then, the root mean square reconstruction error $\mathbf{U}, \{\mathbf{V}_m\}_{m=1}^M$ is 

\begin{equation}
 \text{RMSRE}(\tau) \: =  \: \sqrt{ \Vert \widetilde{\mathbf{A}}_m \Vert_F^2 - \frac{e(\tau)}{M}  \sum_{i=1}^{N} \lambda_i(\mathbf{N})},
\nonumber
\end{equation}

\noindent where $\tau$ is the data compression ratio used in Algorithm \ref{alg1}, $M$ is the number of MC realizations, $N$ is the number of spatial degrees of freedom, $\{\lambda_i(\mathbf{N})\}_{i=1}^N$ denote the eigenvalues of the matrix $\mathbf{N} = \displaystyle \sum_{m=1}^M \widetilde{\mathbf{A}}_m \widetilde{\mathbf{A}}_m^T$ in descending order, and $e(\tau)$ represents the cumulative energy ratio of $\mathbf{N}$:

\begin{equation}
e(\tau) \: := \: \frac{ \sum_{i=1}^{\lceil \tau N \rceil} \lambda_i(\mathbf{N})}{ \sum_{i=1}^N\lambda_i(\mathbf{N})} \: = \: \frac{ \sum_{i=1}^k \lambda_i(\mathbf{N})}{ \sum_{i=1}^N\lambda_i(\mathbf{N})}.
\nonumber
\end{equation}

\end{theorem}

The choice of data compression ratio $\tau$ is crucial for the numerical implementation of our novel low-rank approximation method in Algorithm \ref{alg1}. Selecting a small $\tau$ can effectively reduce the computational and space complexities, while a large $\tau$ generally improves numerical accuracy. Theorem \ref{th3.7} provides a practical strategy for balancing this trade-off. Since the reconstruction error of Algorithm \ref{alg1} decreases as the cumulative energy ratio of $\mathbf{N}$ grows, we can employ the cumulative energy ratio $e(\tau)$ to guide the choice of data compression ratio $\tau$. Actually, $e(\tau)$ serves as a quantitative measure of how much information is retained from the original matrices. From the geometrical perspective, the eigenvectors of a matrix define the principal directions of its corresponding linear transformation, and the eigenvalues describe the scaling magnitude along these directions \cite{mcgivney2014svd}. Consequently, to capture enough information from $\{\widetilde{\mathbf{A}}_m\}_{m=1}^M$, the basis matrix $\mathbf{U}$ should be constructed using a sufficient number of leading eigenvectors of $\mathbf{N}$ in Algorithm \ref{alg1}, so that their associated eigenvalues account for a considerable proportion of the total energy from $\mathbf{N}$. In practice, the optimal data compression ratio $\tau$ is achieved when the cumulative energy ratio $e(\tau)$ is sufficiently close to $1$, which ensures that the low-rank approximations effectively preserve the key features of the original matrices.

Next, we introduce several lemmas that will be used to analyze the error associated with the Neumann series expansion.

\begin{lemma}[\cite{stewart1998matrix}] \label{le3.8}
For a given square matrix $\mathbf{B}$, if its spectral radius satisfies $\rho(\mathbf{B}) < 1$, then the truncation error of the Neumann series is bounded by:

$$\Vert\sum_{r=0}^{R} \mathbf{B}^r - (\mathbf{I}-\mathbf{B})^{-1} \Vert \:\leq \:\frac{\Vert \mathbf{B} \Vert^{R+1}}{1 - \Vert \mathbf{B} \Vert},$$

\noindent for every consistent matrix norm.
\end{lemma}

Lemma \ref{le3.8} also provides a theoretical justification for the perturbation assumption introduced in \eqref{eq3.2}, where $\widetilde{a}$ represents a minor fluctuation around the mean field $\bar{a}$. From a physical perspective, this assumption is physically realistic for porous media applications, where a well-characterized physical process often exhibits small random variations. From a mathematical standpoint, the perturbation framework also improves the tractability of the proposed approach. On the one hand, it effectively guarantees that the total permeability field satisfies the uniform ellipticity condition stated in Assumption \ref{ass4.1}. On the other hand, the small perturbation $\widetilde{a}$ implies that the spectral radius of $\bar{\mathbf{A}}^{-1}\widetilde{\mathbf{A}}_m^*$ remains sufficiently small. According to Lemma \ref{le3.8}, it ensures both the convergence of the Neumann series and the effectiveness of its truncation for matrix inverse approximation.

\begin{lemma}[\cite{wedin1973perturbation}] \label{le3.9}
	  For matrices $\mathbf{A}, \mathbf{B} \in \mathbb{R}^{N \times N}$, if $\text{rank}(\mathbf{A}) = \text{rank}(\mathbf{B})$, then the following bound holds:
	
	\begin{equation} 
		\Vert \mathbf{A}^{-1} - \mathbf{B}^{-1} \Vert \: \le \: \gamma \; \Vert \mathbf{A}^{-1} \Vert_2 \Vert \mathbf{B}^{-1} \Vert_2 \Vert \mathbf{A} - \mathbf{B}\Vert,
		\nonumber
	\end{equation}
	
	\noindent where $\gamma$ is a constant and $\Vert \cdot \Vert$ denotes any consistent matrix norm.
	
\end{lemma}

Finally, the total error estimation of the LRNS solver in Algorithm \ref{alg2} is given as follows.

\begin{theorem}\label{th3.10}
Let $\hat{\bm{u}}_{m,l}$ and $\bm{u}_{m,l}$ be the numerical solution obtained by the LRNS solver in Algorithm \ref{alg2} and the exact solution of \eqref{eq3.4}, respectively. Then, the error estimate between the numerical QoI $\mu_l = \frac{1}{M} \sum_{m=1}^M \hat{\bm{u}}_{m,l}$ and the expectation $\bm{u}_l = \frac{1}{M} \sum_{m=1}^M \bm{u}_{m,l}$ is given by
$$\Vert\,\mu_l \,-\, \bm{u}_l\,\Vert_2 \enspace \leq \enspace C_1 \,\text{RMSRE}(\tau) \:+\: C_2 \,k^{\frac{R+1}{2}},$$

\noindent where $C_1, C_2$ are constants, $\tau$ is the data compression ratio, RMSRE represents the matrix reconstruction error, $k:= \lceil \tau N\rceil$ denotes the reduced dimension, and $R$ is the truncation index for Neumann series.
\end{theorem}

\begin{proof}

By Jensen’s inequality and the properties of norm, the left-hand side reduces to:

\begin{align*}
\Vert\mu_l - \bm{u}_l\Vert_2 \quad & = \quad \Vert \frac{1}{M} \sum_{m=1}^M \hat{\bm{u}}_{m,l} - \frac{1}{M}  \sum_{m=1}^M \bm{u}_{m,l}\Vert_2 \\
& \le \quad \frac{1}{M} \: \sum_{m=1}^M \: \Vert \hat{\bm{u}}_{m,l} - \bm{u}_{m,l} \Vert_2 \\
& = \quad \frac{1}{M} \: \sum_{m=1}^M \: \Vert \left( \bar{\mathbf{A}} + \widetilde{\mathbf{A}}_m\right)^{-1} \bm{b}_l - \sum_{r=0}^R \left(-\bar{\mathbf{A}}^{-1}\mathbf{U}\mathbf{V}_m^T\right)\bar{\mathbf{A}}^{-1}\bm{b}_l\Vert_2\\
& = \quad \frac{1}{M} \: \sum_{m=1}^M \: \Vert \left( \mathbf{I}_N + \bar{\mathbf{A}}^{-1} \widetilde{\mathbf{A}}_m\right)^{-1} \bar{\mathbf{A}}^{-1} \bm{b}_l - \sum_{r=0}^R \left(-\bar{\mathbf{A}}^{-1}\mathbf{U}\mathbf{V}_m^T\right)\bar{\mathbf{A}}^{-1}\bm{b}_l\Vert_2\\
& \le \quad \Vert \bar{\mathbf{A}}^{-1} \Vert_2 \Vert \bm{b}_l \Vert_2 \frac{1}{M} \: \sum_{m=1}^M \: \Vert \left( \mathbf{I}_N + \bar{\mathbf{A}}^{-1} \widetilde{\mathbf{A}}_m\right)^{-1} - \sum_{r=0}^R \left(-\bar{\mathbf{A}}^{-1}\mathbf{U}\mathbf{V}_m^T\right)\Vert_2\\
& \le \quad \Vert \bar{\mathbf{A}}^{-1} \Vert_2 \Vert \bm{b}_l \Vert_2 \frac{1}{M} \: \sum_{m=1}^M \: \Vert \left( \mathbf{I}_N + \bar{\mathbf{A}}^{-1} \widetilde{\mathbf{A}}_m\right)^{-1} - \sum_{r=0}^R \left(-\bar{\mathbf{A}}^{-1}\mathbf{U}\mathbf{V}_m^T\right)\Vert_F\\
& \le \quad \Vert \bar{\mathbf{A}}^{-1} \Vert_2 \Vert \bm{b}_l \Vert_2 \frac{1}{M} \: \sum_{m=1}^M \: \Vert \left( \mathbf{I}_N + \bar{\mathbf{A}}^{-1} \widetilde{\mathbf{A}}_m\right)^{-1} - \left( \mathbf{I}_N + \bar{\mathbf{A}}^{-1} \mathbf{U}\mathbf{V}_m^T\right)^{-1} \\
& \qquad \qquad + \left( \mathbf{I}_N + \bar{\mathbf{A}}^{-1} \mathbf{U}\mathbf{V}_m^T\right)^{-1} - \sum_{r=0}^R \left(-\bar{\mathbf{A}}^{-1}\mathbf{U}\mathbf{V}_m^T\right)\Vert_F\\
& \le \quad \Vert \bar{\mathbf{A}}^{-1} \Vert_2 \Vert \bm{b}_l \Vert_2 \frac{1}{M} \: \sum_{m=1}^M \: \Vert \left( \mathbf{I}_N + \bar{\mathbf{A}}^{-1} \widetilde{\mathbf{A}}_m\right)^{-1} - \left( \mathbf{I}_N + \bar{\mathbf{A}}^{-1} \mathbf{U}\mathbf{V}_m^T\right)^{-1}\Vert_F \\
& \qquad + \Vert \bar{\mathbf{A}}^{-1} \Vert_2 \Vert \bm{b}_l \Vert_2 \frac{1}{M} \: \sum_{m=1}^M \: \Vert \left( \mathbf{I}_N + \bar{\mathbf{A}}^{-1} \mathbf{U}\mathbf{V}_m^T\right)^{-1} - \sum_{r=0}^R \left(-\bar{\mathbf{A}}^{-1}\mathbf{U}\mathbf{V}_m^T\right)\Vert_F.
\end{align*}

According to Lemma \ref{le3.9}, we can derive the following formulation

\begin{align*}
&\frac{1}{M^2} \: \sum_{m=1}^M \: \Vert \left( \mathbf{I}_N + \bar{\mathbf{A}}^{-1} \widetilde{\mathbf{A}}_m\right)^{-1} - \left( \mathbf{I}_N + \bar{\mathbf{A}}^{-1} \mathbf{U}\mathbf{V}_m^T\right)^{-1}\Vert_F^2\\
\le\quad&\frac{\gamma^2}{M^2} \: \sum_{m=1}^M \: \Vert \left( \mathbf{I}_N + \bar{\mathbf{A}}^{-1} \widetilde{\mathbf{A}}_m\right)^{-1} \Vert_2^2 \;\Vert\left( \mathbf{I}_N + \bar{\mathbf{A}}^{-1} \mathbf{U}\mathbf{V}_m^T\right)^{-1}\Vert_2^2 \; \Vert\bar{\mathbf{A}}^{-1} \widetilde{\mathbf{A}}_m - \bar{\mathbf{A}}^{-1} \mathbf{U}\mathbf{V}_m^T\Vert_F^2\\
\le\quad&\frac{\gamma^2}{M^2} \: \sum_{m=1}^M \: \Vert \left( \mathbf{I}_N + \bar{\mathbf{A}}^{-1} \widetilde{\mathbf{A}}_m\right)^{-1} \Vert_2^2 \;\Vert\left( \mathbf{I}_N + \bar{\mathbf{A}}^{-1} \mathbf{U}\mathbf{V}_m^T\right)^{-1}\Vert_2^2 \; \Vert \bar{\mathbf{A}}^{-1} \Vert_F^2 \;\Vert \widetilde{\mathbf{A}}_m - \mathbf{U}\mathbf{V}_m^T\Vert_F^2\\
\le\quad&\frac{\gamma^2}{M} \: \sum_{m=1}^M \: \Vert \left( \mathbf{I}_N + \bar{\mathbf{A}}^{-1} \widetilde{\mathbf{A}}_m\right)^{-1} \Vert_2^2 \;\Vert\left( \mathbf{I}_N + \bar{\mathbf{A}}^{-1} \mathbf{U}\mathbf{V}_m^T\right)^{-1}\Vert_2^2 \; \Vert \bar{\mathbf{A}}^{-1} \Vert_F^2 \; \text{RMSRE}(\tau)^2.\\
\end{align*}

Therefore, the first term of the total error estimate above reduces to

\begin{equation}
\Vert \bar{\mathbf{A}}^{-1} \Vert_2 \Vert \bm{b}_l \Vert_2 \frac{1}{M} \: \sum_{m=1}^M \: \Vert \left( \mathbf{I}_N + \bar{\mathbf{A}}^{-1} \widetilde{\mathbf{A}}_m\right)^{-1} - \left( \mathbf{I}_N + \bar{\mathbf{A}}^{-1} \mathbf{U}\mathbf{V}_m^T\right)^{-1}\Vert_F \:
\le\: C_1 \; \text{RMSRE}(\tau),\\
\nonumber
\end{equation}

\noindent where $C_1:=\frac{\gamma}{\sqrt{M}}\Vert \bar{\mathbf{A}}^{-1} \Vert_F\Vert \bar{\mathbf{A}}^{-1} \Vert_2 \Vert \bm{b}_l \Vert_2 \mathop{\max}\limits_{m = 1,\cdots,M}\left\{\Vert \left( \mathbf{I}_N + \bar{\mathbf{A}}^{-1} \widetilde{\mathbf{A}}_m\right)^{-1} \Vert_2^2 \;\Vert\left( \mathbf{I}_N + \bar{\mathbf{A}}^{-1} \mathbf{U}\mathbf{V}_m^T\right)^{-1}\Vert_2^2\right\}$ is a constant. By Theorem \ref{th3.7}, the RMSRE value arising from Algorithm \ref{alg1} depends
on the data compression ratio $\tau$ and it is proportional to the square root of the cumulative energy ratio $e(\tau)$ of the matrix $\mathbf{N}$.

Based on Theorem \ref{th3.2}, Lemma \ref{le3.8} and the properties of matrix norm, the second term of the total error estimate is bounded by

\begin{align*}
& \Vert \bar{\mathbf{A}}^{-1} \Vert_2 \Vert \bm{b}_l \Vert_2 \frac{1}{M} \: \sum_{m=1}^M \: \Vert \left( \mathbf{I}_N + \bar{\mathbf{A}}^{-1} \mathbf{U}\mathbf{V}_m^T\right)^{-1} - \sum_{r=0}^R \left(-\bar{\mathbf{A}}^{-1}\mathbf{U}\mathbf{V}_m^T\right)\Vert_F\\
\le \quad &\Vert \bar{\mathbf{A}}^{-1} \Vert_2 \Vert \bm{b}_l \Vert_2 \frac{1}{M} \: \sum_{m=1}^M \;\frac{\Vert -\bar{\mathbf{A}}^{-1}\mathbf{U}\mathbf{V}_m^T \Vert_F^{R+1}}{1-\Vert -\bar{\mathbf{A}}^{-1}\mathbf{U}\mathbf{V}_m^T \Vert_F}\\
\le \quad &\Vert \bar{\mathbf{A}}^{-1} \Vert_2 \Vert \bm{b}_l \Vert_2 \frac{1}{M} \: \sum_{m=1}^M \;\Vert -\bar{\mathbf{A}}^{-1}\mathbf{U}\mathbf{V}_m^T \Vert_F^{R+1}\\
= \quad &\Vert \bar{\mathbf{A}}^{-1} \Vert_2 \Vert \bm{b}_l \Vert_2 \frac{1}{M} \: \sum_{m=1}^M \;\Vert -\bar{\mathbf{A}}^{-1}\mathbf{U}\mathbf{U}^T \widetilde{\mathbf{A}}_m\Vert_F^{R+1}\\
\le \quad & \frac{1}{M} \Vert \bar{\mathbf{A}}^{-1} \Vert_2 \Vert \bm{b}_l \Vert_2 \Vert \bar{\mathbf{A}}^{-1}\Vert_F^{R+1} \: \sum_{m=1}^M \;\Vert \widetilde{\mathbf{A}}_m \Vert_F^{R+1}\;\Vert \mathbf{U}\mathbf{U}^T\Vert_F^{R+1}.
\end{align*}

Denote the constant $C2 := \frac{\Vert \bm{b}_l \Vert_2}{M} \Vert \bar{\mathbf{A}}^{-1} \Vert_2  \Vert \bar{\mathbf{A}}^{-1}\Vert_F^{R+1} \: \sum_{m=1}^M \;\Vert \widetilde{\mathbf{A}}_m \Vert_F^{R+1}$. Then we can further simplify this term based on the definition of Frobenius norm and the orthogonality of $\mathbf{U}$:

\begin{align*}
& \Vert \bar{\mathbf{A}}^{-1} \Vert_2 \Vert \bm{b}_l \Vert_2 \frac{1}{M} \: \sum_{m=1}^M \: \Vert \left( \mathbf{I}_N + \bar{\mathbf{A}}^{-1} \mathbf{U}\mathbf{V}_m^T\right)^{-1} - \sum_{r=0}^R \left(-\bar{\mathbf{A}}^{-1}\mathbf{U}\mathbf{V}_m^T\right)\Vert_F \quad
\le \quad  C_2\;\Vert \mathbf{U}\mathbf{U}^T\Vert_F^{R+1}\\
& \qquad \quad \le \enspace C_2\; \text{tr}\left( (\mathbf{U}\mathbf{U}^T)^T (\mathbf{U}\mathbf{U}^T)\right)^{\frac{R+1}{2}} \enspace = \enspace C_2\; \text{tr}\left( \mathbf{U}^T\mathbf{U}\right)^{\frac{R+1}{2}}\enspace = \enspace C_2\; \text{tr}\left( \mathbf{I}_k\right)^{\frac{R+1}{2}} \enspace = \enspace C_2\; k^{\frac{R+1}{2}},\\
\end{align*}

\noindent where $k := \lceil \tau N \rceil$ is the reduced dimension used in Algorithm \ref{alg1}.
\end{proof}

\section{Application to unsteady stochastic diffusion PDEs and the unsteady stochastic diffusion optimal control problems} \label{sec4}

In this section, we apply the LRNS solver in Algorithm \ref{alg2} to two specific UQ problems: the unsteady diffusion equations with a random permeability coefficient and the optimal control problem governed by the unsteady stochastic diffusion equations.

\subsection{Unsteady stochastic diffusion partial differential equations} \label{subsec4.1}

Let $\mathcal{D} \subset \mathbf{R}^d$, $d=1,2,3$, be a bounded, Lipschitz domain and let $\mathcal{Q}$ denote a space-time cylinder with $T>0$. We consider the following unsteady diffusion-type SPDE: find a stochastic function $u : \mathcal{D} \times [0,T] \times \Omega \rightarrow \mathbf{R}$ such that

\begin{equation}\label{eq4.1}
\left\{\:
\begin{aligned}
\partial_t u(\bm{x}, t, \omega) - \nabla \cdot \big( a(\bm{x},\omega) \nabla u(\bm{x}, t, \omega) \big) &= f(\bm{x}, t), \quad \quad\text{in} \enspace \mathcal{D} \times [0,T] \times \Omega,\\
u(\bm{x}, t, \omega) &=g(\bm{x},t),  \qquad \text{on}\enspace \partial \mathcal{D} \times [0,T] \times \Omega,\\
u(\bm{x}, 0, \omega) &= u_0(\bm{x}),  \qquad \,\,\text{on} \enspace \mathcal{D} \times \{t=0\}\times\Omega,
\end{aligned}
\right.
\end{equation}

\noindent where the diffusion coefficient $a$ is an almost surely continuous and positive random field on $\mathcal{D}$. We assume for simplicity that the deterministic force term, boundary function and initial condition $f,g,u_0 \in \mathcal{L}^2(\mathcal{D})$.

To guarantee the existence and uniqueness of the solution of \eqref{eq4.1}, we impose the following condition on the stochastic diffusion coefficient $a(\cdot, \omega)$.
	
\begin{assumption} [Uniform ellipticity condition] \label{ass4.1}

Let $a_{min} := ess \inf_{\omega} \Vert a(\cdot, \omega) \Vert_{\mathcal{H}_0^1(\mathcal{D})}$ and $a_{max} := ess \sup_{\omega} \Vert a(\cdot, \omega) \Vert_{\mathcal{H}_0^1(\mathcal{D})}$. Then, the stochastic diffusion coefficient $a$ satisfies that
		
\begin{equation}
	0 \, < \, a_{\min} \leq a(\bm{x}, \omega) \leq a_{\max} \, < \, \infty, \quad \text{ for a.e. }(\bm{x}, \omega)  \in \mathcal{D} \times \Omega.
	\nonumber
\end{equation}
		
\end{assumption} 

In practical engineering and geoscientific applications, permeability is typically inferred using data-driven approaches, such as regression analysis and Bayesian learning techniques \cite{roding2020predicting}. Motivated by the practical modeling considerations, we adopt the additive perturbation form in \eqref{eq3.2} and represent the diffusion coefficient using the truncated Karhunen-Loe\.{v}e (KL) expansion \cite{ghanem2003stochastic}:

\begin{equation}\label{eq4.3} 
a(\bm{x},\omega) \:=\: \bar{a}(\bm{x}) \:+\: \sigma\sum_{t=1}^{T} \sqrt{\lambda_t}\, r_t(x) Y_t(\omega),
\end{equation}

\noindent where $\bar{a}(\bm{x})$ denotes the expectation of the random field $a(\bm{x},\omega)$:

$$\bar{a}(\bm{x}) \:=\: \mathbb{E}\left[ a(\bm{x},\omega)\right] \::=\: \int_\Omega a(\bm{x},\omega) \mathrm{d}\mathbb{P}(\omega).$$

\noindent $\{\lambda_t\}_{t=1}^T$, $\{r_t(x)\}_{t=1}^\infty$ are the positive eigenvalues in decreasing order, and the corresponding orthonormal eigen-functions of the covariance function:

$$\text{Cov}(\bm{x},\bm{x}^{\prime}) \:=\: \int_\Omega \left(a(\bm{x},\omega) - \bar{a}(\bm{x})\right) \left(a(\bm{x}^{\prime},\omega) - \bar{a}(\bm{x}^{\prime})\right)\mathrm{d}\mathbb{P}(\omega),$$

\noindent respectively. The random variables $\{Y_t(\omega)\}_{t=1}^{T}$ in \eqref{eq4.3} are assumed to be mutually independent with zero mean and unit variance. The scaling index $\sigma \in [0,1]$ controls the magnitude of the perturbation $\widetilde{a}(\bm{x},\omega)$, which ensures that the resulting permeability field satisfies the uniform ellipticity condition.

Conventionally, the random variables $\{Y_t(\omega)\}_{t=1}^T$ in the truncated KL expansion follow the standard normal distribution. However, a finite sum of Gaussian random variables may produce negative or unbounded permeability values with small probability. To satisfy the uniform ellipticity condition, we model the underlying KL coefficients $\{Y_t(\omega)\}_{t=1}^T$ as a sequence of independent and identically distributed random variables following the \textbf{truncated standard normal distribution} \cite{burkardt2014truncated, gunzburger2016stochastic}, which ensures the boundedness of the perturbation $\widetilde{a}(\bm{x},\omega)$. Although it may still take negative values, the dominant mean field $\bar{a}>0$ and the scaling index $\sigma$ guarantee that the total permeability $a = \bar{a} + \widetilde{a}$ remains strictly positive for all realizations. By Chebyshev's inequality, approximately $99.73\%$ of the values of a normal distribution lie within three standard deviations of its mean. Therefore, the truncated standard normal distribution preserves essential statistical features of the standard normal distribution, while avoiding extreme and non-physical values. 

Finally, we define weak solutions of \eqref{eq4.1}: find $u \in \mathcal{H}^1(\mathcal{Q})$ such that for a.e. $0 \leq t \leq T$,

\begin{equation}\label{eq4.2}
\mathbb{E}\left[ (u_t,v) \:+\: a(u,v)\right]\enspace=\enspace (f,v), \quad \forall v \in \mathcal{H}_0^1(\mathcal{Q}),
\end{equation}

\noindent where the bilinear forms are defined as 
$$a(u,v) \enspace = \enspace \int_\mathcal{D} a \nabla u \cdot \nabla v \,\mathrm{d}\bm{x},$$
\noindent and
$$(u,v) \enspace = \enspace \int_\mathcal{D} \,u \, v \,\mathrm{d}\bm{x}.$$

The following lemma presents the existence and uniqueness of the solution of the weak form \eqref{eq4.2}.
	
\begin{lemma} [\cite{evans2022partial}] \label{le4.2}		
Let the source $f \in \mathcal{L}^2(\mathcal{Q})$ and Assumption \ref{ass4.1} hold. There exists a unique solution $u$ of the weak formulation \eqref{eq4.2} in the Hilbert space $\mathcal{H}^2(\mathcal{Q})$.
\end{lemma}

\subsubsection{Discretization of the weak formulation} \label{subsubsec4.1.2}

We begin by using the standard MC-FEM method to discretize both the physical and probabilistic spaces. This approach is effective in mitigating the curse of dimensionality. In practice, the expectation operator $\mathbb{E}$ is approximated by sample averages of i.i.d. realizations, while the spatial domain is discretized using the finite element method.

Consider a regular mesh $\mathcal{T}_h$ over $\mathcal{D}$ with the mesh size $h>0$. Let $V_h \subset \mathcal{H}^1(\mathcal{D})$ be the finite element space with the basis functions $\{\phi_j\}_{j=1}^N$. Then we can express the approximate solution as follows:

\begin{equation} 
	u_h (\bm{x},t,\omega) = \sum_{j=1}^N u_j(t,\omega) \phi_j(\bm{x}).
\nonumber
\end{equation}

We generate $M$ MC samples $\widetilde{a}_m$ of the perturbation field in the permeability coefficient. Thus, the semi-discrete system resulting from the standard MC-FEM approach is given by

\begin{equation} 
\begin{aligned}
\sum_{j=1}^N \frac{\mathrm{d}u_{j,m}}{\mathrm{d}t} \int_D \phi_j \phi_i d\bm{x} \:&+\:\sum_{j=1}^N u_{j,m} \left[\int_D \bar{a} \nabla \phi_j \nabla \phi_i d\bm{x}  + \:  \int_D \widetilde{a}_m \nabla \phi_j \nabla \phi_i d\bm{x} \right]\\
& \: = \: \int_D f \phi_i d\bm{x}, \quad i = 1,\cdots,N, m = 1,\cdots,M.
\nonumber
\end{aligned}
\end{equation}

Equivalently, it constitutes a linear system of algebraic equations:

\begin{equation} \label{eq4.5}
\mathbf{G} \frac{\mathrm{d}}{\mathrm{d}t}\bm{u}_m(t) \:+\: \left(\bar{\mathbf{A}} + \widetilde{\mathbf{A}}_m\right)\bm{u}_m(t) \:=\: \bm{b}(t), \quad t \in (0,T], 
\end{equation}

\noindent for $m = 1,\cdots,M$. The mass matrix $\mathbf{G} \in \mathbb{R}^{N \times N}$, the stiffness matrices $\bar{\mathbf{A}}, \widetilde{\mathbf{A}}_m \in \mathbb{R}^{N \times N}$ and the load vector $\bm{b} \in \mathbb{R}^{N \times 1}$ are respectively given by

\begin{equation} 
\begin{aligned}\label{eq4.6}
&\mathbf{G}_{ij} = \int_\mathcal{D} \phi_j \phi_i d\bm{x}, \quad \bar{\mathbf{A}}_{ij} = \int_D \bar{a} \nabla \phi_j \nabla \phi_i d\bm{x},\quad \left(\widetilde{\mathbf{A}}_m\right)_{ij} = \int_D \widetilde{a}_m \nabla \phi_j \nabla \phi_i d\bm{x}, \\
&\bm{b}_i(t) = \int_D f(t) \phi_i d\bm{x},\quad\bm{u}_m(t) = \left( u_{1,m}(t), ..., u_{N,m}(t)\right) ^T, \qquad m = 1,\cdots,M.
\end{aligned}
\end{equation}

For temporal discretization, we partition $(0,T]$ into $L$ sub-intervals with the time step $\Delta t = T/L$ and time nodes $t_l = l\Delta t$, $l=0,\dots,L$. The Crank-Nicolson scheme is applied to the linear system \eqref{eq4.5} and we obtain the following full discrete system: for each sample $m = 1,\cdots,M$ and time step $l = 1,\cdots,L$, find $u_h^{m,n} \in V_h$ such that

\begin{equation} 
\mathbf{G} \frac{\bm{u}_{m,l+1} - u_{m,l}}{\Delta t} + \frac{1}{2}\left(\bar{\mathbf{A}} + \widetilde{\mathbf{A}}_m\right)\bm{u}_{m,l+1} + \frac{1}{2}\left(\bar{\mathbf{A}} + \widetilde{\mathbf{A}}_m\right)\bm{u}_{m,l} \:=\: \frac{1}{2}\bm{b}_{l+1} + \frac{1}{2}\bm{b}_{l}, 
\nonumber
\end{equation}

\noindent and it further reduces to

\begin{equation} \label{eq4.7}
\left(\bar{\mathbf{K}} + \widetilde{\mathbf{K}}_m\right) \bm{u}_{m,l+\frac{1}{2}} \enspace=\enspace \bm{f}_{m,l+\frac{1}{2}},
\end{equation}

\noindent with the initial state $\bm{u}(0,\omega_m) = \bm{u}_0$. Then, the matrices and vector are defined as below:

\begin{equation}\label{eq4.8}
\begin{aligned}
&\bar{\mathbf{K}} \::=\: \frac{\mathbf{G}}{\Delta t/2} + \bar{\mathbf{A}}, \quad \widetilde{\mathbf{K}}_m \::=\: \widetilde{\mathbf{A}}_m,\quad \bm{f}_{m,l+\frac{1}{2}} \::=\:\frac{1}{2}\bm{b}_{l+1} + \frac{1}{2}\bm{b}_{l} + \frac{\mathbf{G}}{\Delta t/2}\bm{u}_{m,l},\\
& \bm{u}_{m,l+1} \: =\: 2\bm{u}_{m,l+\frac{1}{2}}-\bm{u}_{m,l},\qquad m = 1,\cdots,M, \: l = 0,\cdots,L.
\end{aligned}
\end{equation}

\subsubsection{Construct the LRNS solver for the unsteady stochastic diffusion PDE} \label{subsubsec4.1.3}

Solving the algebraic equations in \eqref{eq4.8} for all the MC realizations and time steps leads to heavy computational and memory burdens. Therefore, to address these issues, our LRNS solver first finds the low-rank representation of the perturbation matrices $\{\widetilde{\mathbf{K}}_m\}_{m=1}^M$ using Algorithm \ref{alg1}:

\begin{equation}\label{eq4.9}
\widetilde{\mathbf{K}}_m \enspace\approx\enspace \mathbf{U} \mathbf{V}_m^T,\quad m = 1,\cdots,M.
\end{equation}

\noindent Then, we utilize the Neumann matrix series to approximate the inverse of coefficient matrices, and the numerical sample solutions can be expressed as

\begin{equation}\label{eq4.10}
\left\{\:
\begin{aligned}
& \mathbf{u}_{m,l+\frac{1}{2}} = \sum_{r=0}^R \left(-\bar{\mathbf{K}}^{-1}\mathbf{U}\mathbf{V}_m^T\right)^{r} \bar{\mathbf{K}}^{-1}\bm{f}_{m,l+\frac{1}{2}}, \\
& \bm{u}_{m,l+1} \: =\: 2\bm{u}_{m,l+\frac{1}{2}}-\bm{u}_{m,l},
\end{aligned}
\right.
\end{equation}

\noindent for each $m=1,2,\cdots,M$ and $l = 1,2,\cdots, L$. $R\ge0$ is a given truncation index for the Neumann series. Finally, the quantity of interest to approximate the expected solution at time $t_l$ is given by

\begin{equation} \label{eq4.11}
\mu(\bm{x}, t_l) \enspace:=\enspace  \frac{1}{M} \sum_{m=1}^M \sum_{j=1}^N u_{j,m,l} \,\phi_j(\bm{x}).
\end{equation}

Therefore, we derive the efficient LRNS solver for the unsteady diffusion PDE with a random permeability coefficient. Instead of repeatedly inverting the $N$-dimensional matrices, the formulation in \eqref{eq4.9} reduces the problem to a sequence of matrix–vector multiplications, thereby effectively accelerating the overall computations. In addition, the low-rank approximation in \eqref{eq4.8} significantly reduces storage requirements for the large-scale stiffness matrices. The complete procedure of the LRNS solver for the unsteady stochastic diffusion equation is summarized in Algorithm \ref{alg3}.

\begin{algorithm}[h]
\caption{The LRNS solver for the unsteady stochastic diffusion equation with random permeability in \eqref{eq4.1}}
\label{alg3}
\begin{algorithmic}[1]
\Require A tessellation $\mathcal{T}_h$ of spatial domain $\mathcal{D}$, final time $T$, MC sample size $M$, numbers of spatial and temporal discretization for SPDEs $N$ and $L$, Neumann series truncation index $R$, and data compression ratio $\tau$.
\Ensure Approximation of QoI $\mu$.

\State Construct the finite element space $V_h \subset \mathcal{H}^1(\mathcal{D})$.

\State Generate the deterministic and perturbation permeability $\bar{a},\{\widetilde{a}_m\}_{m = 1}^M$ by the truncated KL expansion in \eqref{eq4.3}.
		
\State Assemble the stiffness matrices $\bar{\mathbf{K}}, \{\mathbf{\widetilde{K}}_m\}_{m=1}^M$ and the load vectors $\bm{f}_{m,l+\frac{1}{2}}$ by \eqref{eq4.8}.
	
\State Deal with the Dirichlet boundary conditions of $\bar{\mathbf{K}}$.

\noindent \hspace{-2em} \textit{Step 2: Low-rank approximation and matrix series inversion}

\State Determine the reduced dimension $k = \lceil \tau N \rceil$ and compute the optimal rank-$k$ approximation $\mathbf{U}, \{\mathbf{V}_m\}_{m=1}^M$ of $\{\widetilde{\mathbf{K}}_m\}_{m=1}^M$ using Algorithm \ref{alg1}.
 
\For{$l = 0, \cdots, L$}
\For{$m = 1, \ldots, M$} 
    \State Compute the sample solution $\bm{u}_{m,l+1}$ by the low-rank correction form in \eqref{eq4.10}.
\EndFor
\EndFor

\State \Return the QoI $\mu$ by \eqref{eq4.11}.
\end{algorithmic}
\end{algorithm}

\subsubsection{Numerical experiments on the unsteady stochastic diffusion equation} \label{subsubsec4.1.4}

In this subsection, we conduct numerical experiments on the unsteady stochastic diffusion PDE in \eqref{eq4.1} to evaluate the numerical performance of the proposed LRNS solver. The simulations are carried out in MATLAB R2022a software on an Apple M1 machine with 8 GB of memory. 

We consider a two-dimensional rectangle domain $\mathcal{D} = [0,1]\times[0,1]$ over the time interval $t \in [0,1]$. The numerical example is tested under the following settings. Let the spatial variable $\bm{x} = [x,y]^T$. The source term, boundary condition function, and initial condition function in \eqref{eq4.1} are respectively given by

\begin{equation}
f(\bm{x},t) = 1, \quad g(\bm{x},t) = 0,\quad u_0(\bm{x}) = \sin(2\pi x)\sin(2\pi y).
\nonumber
\end{equation}

For the spatial discretization, we employ a uniform tessellation $\mathcal{T}_h$ with the mesh size $h = 1/16$ and the degrees of freedom $N = 1089$. Numerical integration is performed using the nine-point Gauss quadrature rule in the experiments. The MC sampling method is adopted to discretize the probabilistic space with the sample size $M = 1000$. The temporal domain is discretized using a time step $\Delta t = 0.01$ and the number of time sub-intervals $L = 100$. We adopt the truncation index for the Neumann matrix series $R=5$ in the numerical implementation of the proposed algorithms. the truncated KL expansion in \eqref{eq4.3}. The mean permeability, the scaling index, the parameter $T$, and the covariance function for the truncated KL expansion in \eqref{eq4.3} are defined by:
$$\bar{a}(\bm{x}) = 1, \:\sigma = 0.2,\:T=19,\:\text{Cov}(\bm{x},\bm{y})  \:= \:  \exp\left( -\frac{\vert \bm{x}-\bm{y}\vert}{0.2}\right).$$ 
\noindent The random variables $\{Y_t(\omega)\}_{t=1}^T$ follow the standard normal distribution $\mathcal{N}(0,1)$ with the truncated interval $\left[ -3,3\right] $. According to the three-sigma principle, this truncation retains the main statistical features of the original distribution and keeps the permeability uniformly bounded, thereby ensuring the uniform ellipticity condition.

\begin{figure}[!htbp]
\centering
\includegraphics[width=1\textwidth]{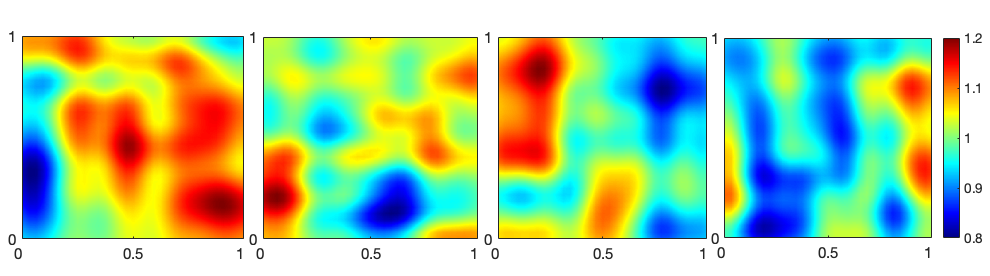}
\caption{Four randomly selected MC realizations of the permeability $a(\bm{x},\omega)$.}
\label{fig4.2}
\end{figure}

\begin{figure}[htbp]
\centering
\includegraphics[width=1\textwidth]{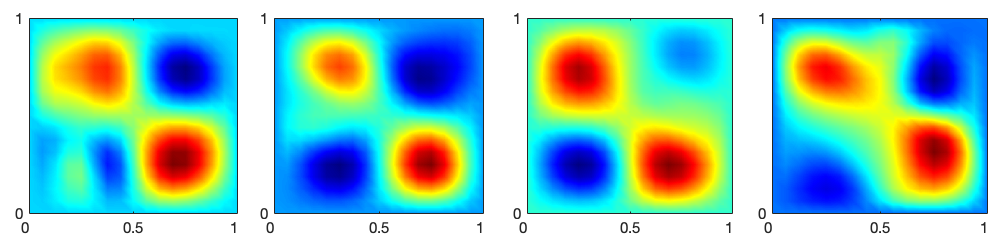}
\caption{Four randomly selected samples of reference solution.}
\label{fig4.4}
\end{figure}

Four randomly selected realizations of $a(\bm{x},\omega)$ are displayed in Figure \ref{fig4.2}, which exhibits the randomness of the diffusion coefficient. The reference solution is generated through the standard MC-FEM method with the reference tic-toc time $t_{ref} = 2887.097$. Figure \ref{fig4.4} demonstrates the sample solutions with four randomly selected MC realizations of $a(\bm{x},\omega)$, which manifests the perturbative nature of the system.

\begin{figure}[htbp] 
\centering
\includegraphics[width=1\textwidth]{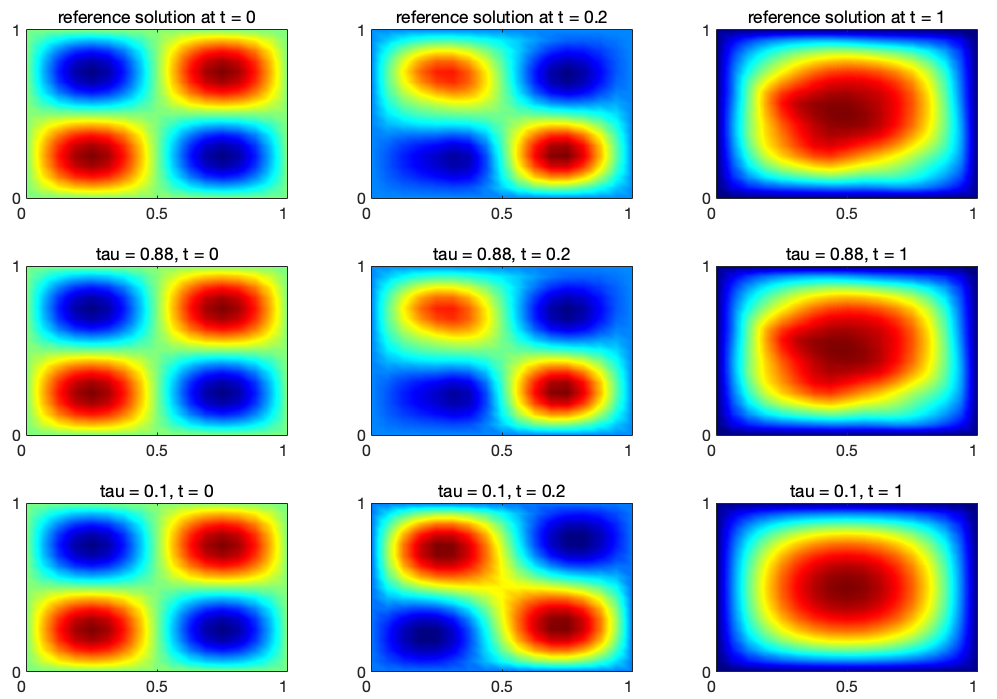}
\caption{Comparison among the reference solution, the numerical solutions computed by the LRNS solver in Algorithm \ref{alg3} under data compression ratios $\tau = 88\%$ and $10\%$.}
\label{fig4.5}
\end{figure}

Then we assess the numerical accuracy and efficiency of the LRNS solver in Algorithm \ref{alg3} under different choices of data compression ratios $\tau = 100\%, 95\%, 88\%, 50\%, 30\%, 10\%$. Figure \ref{fig4.5} and Table \ref{tab4.1} demonstrate that when adopting an appropriate choice of data compression ratio, $88\%$ in this experiment, our LRNS solver can maintain a high numerical precision, and reduce computational and space complexities associated with large-scale matrices at the same time, which confirms the effectiveness and validity of the proposed algorithm.

\begin{table}[htbp]
\centering 
\caption{Simulation results for Algorithm \ref{alg3} about the reduced dimensions, the numerical error, and the tic-toc time under six different data compression ratios $\tau = 100\%, 95\%, 88\%, 50\%, 30\%, 10\%$.}
\resizebox{\textwidth}{!}{
\begin{tabular}{cccccccc}
\hline
Data compression ratio $\tau$ & $100\%$ & $95\%$ & $88\%$ & $50\%$ & $30\%$ & $10\%$ \\\hline
Reduced dimension $k$ & $1089$ & $1035$ & $959$ & $545$ & $327$ & $109$\\
Error & $6.5169 \times 10^{-5}$ & $6.7488 \times 10^{-5}$ & $6.8119 \times 10^{-5}$ & $3.1163 \times 10^{-3}$ &  $7.3126 \times 10^{-2}$ & $1.1329 \times 10^{-1}$ \\
Tic-toc time (s) & $2769.970$ & $2619.275$ & $2301.890$ & $2178.805$ & $2051.835$ & $1722.615$\\
$t_{num}/t_{ref}$ & $95.94\%$ & $90.72\%$ & $79.73\%$ & $75.47\%$ & $71.07\%$ & $59.67\%$\\\hline 
\end{tabular}}
\label{tab4.1}
\end{table}

\begin{figure}[htbp]
\centering
\includegraphics[width=1\textwidth]{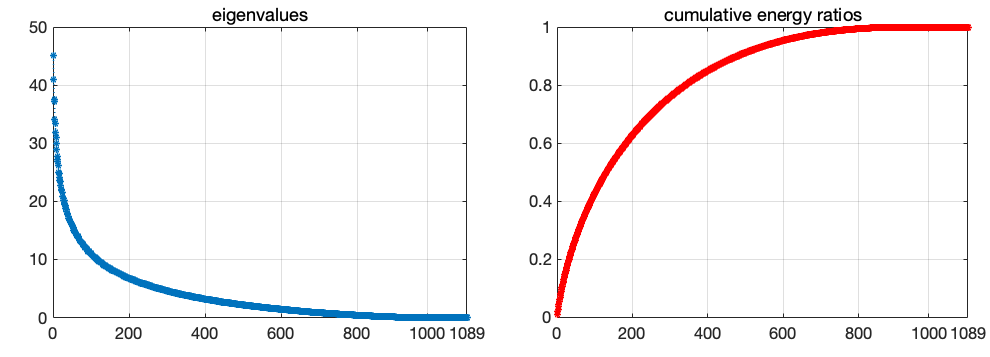}
\caption{The eigenvalues of the matrix $\mathbf{N} = \sum_{m=1}^M \widetilde{\mathbf{A}}_m \widetilde{\mathbf{A}}_m^T$ in descending order (left) and their corresponding cumulative energy ratios (right).}
\label{fig4.6}
\end{figure}

As stated in Section \ref{sec3}, a large data compression ratio generally enhances the accuracy of the LRNS solver, while a small $\tau$ leads to significant savings in memory and run-time. Theorem \ref{th3.7} provides practical guidance for determining a proper data compression ratio based on the cumulative energy ratio of the matrix $\mathbf{N}$. According to the simulation results in Figure \ref{fig4.6}, we select the data compression ratio $\tau = 87.97\%$ and the reduced dimension $k = 958$, when the corresponding cumulative energy ratio $e(87.97\%)$ is extremely close to $1$.

\begin{figure}[htbp]
\centering
\includegraphics[width=0.7\textwidth]{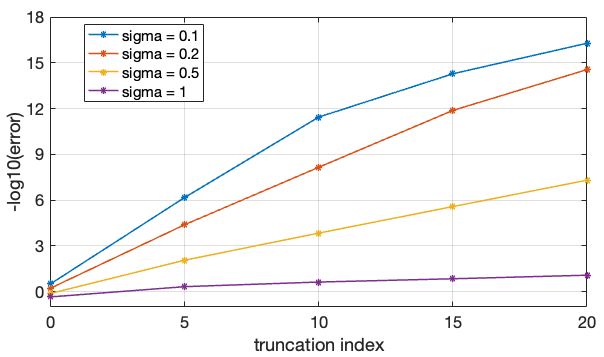}
\caption{The log-error values of the LRNS solver under the truncation index for the Neumann series $R=0,5,10,15$ and the scaling index $\sigma = 0.1,0.2,0.5,1$}
\label{fig4.7}
\end{figure}

Finally, we examine the sensitivity of the LRNS solver to the scaling index. Figure \ref{fig4.7} compares the numerical errors obtained by the LRNS solver for four different $\sigma = 0.1,0.2,0.5,1$. For the smallest choice, $\sigma = 0.1$, Algorithm \ref{alg3} achieves the accuracy on the order of $10^{-7}$ using only $R=5$ terms in the Neumann matrix series. However, as $\sigma$ increases, a larger truncation index $R$ is required to attain comparable precision. This observation is consistent with the conclusion in Lemma \ref{le3.8}. When the model exhibits strong random variability, the spectral radius of $\bar{\mathbf{A}}^{-1}\widetilde{\mathbf{A}}_m^*$ increases, which leads to a slow convergence rate of the associated Neumann series.

\subsection{Unsteady stochastic diffusion optimal control problems} \label{subsec4.2}

In this section, we extend the LRNS solver to the stochastic optimal control framework. We investigate the minimization of the following space-time tracking cost functional:

\begin{equation} \label{eq4.12}
\min \mathcal{J}(u, f) \: := \: \bm{E} \left[\frac{1}{2} \int_0^T \Vert u - U \Vert^2 \mathrm{d}t \right] + \frac{\beta}{2} \int_0^T \Vert f \Vert^2 \mathrm{d}t,
\end{equation}

\noindent with respect to the stochastic state $u$ and the deterministic control $f$ subject to the stochastic parabolic initial boundary value problem:
	
\begin{equation}\label{eq4.13}
\left\{\:
\begin{aligned}
\partial_t u(\bm{x}, t, \omega) - \nabla \cdot \big( a(\bm{x},\omega) \nabla u(\bm{x}, t, \omega) \big) &= f(\bm{x}, t), \quad \quad\text{in} \enspace \mathcal{D} \times [0,T] \times \Omega,\\
u(\bm{x}, t, \omega) &=g(\bm{x},t),  \qquad \text{on}\enspace \partial \mathcal{D} \times [0,T] \times \Omega,\\
u(\bm{x}, 0, \omega) &= u_0(\bm{x}),  \qquad \,\,\text{on} \enspace \mathcal{D} \times \{t=0\}\times\Omega,
\end{aligned}
\right.
\end{equation}
	
\noindent where $U \in \mathcal{L}^2(\mathcal{Q})$ is the given deterministic desired state and $\beta>0$ is a small control penalty parameter. The admissible set $\mathcal{U}_{ad}$ is defined as
	
\begin{equation}
\mathcal{U}_{ad} \enspace := \enspace L^2(0,T; H^2(\mathcal{D}) \cap H_0^1(\mathcal{D}))\: \otimes \: \{ f \in L^2(D): f(\bm{x}) \ge 0, \forall \bm{x} \in \mathcal{D} \}.
\nonumber
\end{equation}
	
In this setting, the standard weak formulation of the distributed optimal control problem in \eqref{eq4.12} and \eqref{eq4.13} reads as follows: find the optimal state $u^*$ and control $f^*$ satisfying

\begin{equation} \label{eq4.14}
\mathop{\min}\limits_{(u,f) \: \in \: \mathcal{U}_{ad}}  \: \mathcal{J}(u,f), \qquad s.t. \enspace \mathbb{E}\left[ (u_t,v) \:+\: a(u,v)\right]\:=\: (f,v), \quad \forall v \in \mathcal{H}_0^1(\mathcal{Q}).
\end{equation}

The well-posedness of the weak solution of \eqref{eq4.14} is present in the following lemma.
	
\begin{lemma} [\cite{gunzburger2011space}] \label{le4.3}
If the admissible set $\mathcal{U}_{ad}$ is nonempty, then the distributed optimal control problem in \eqref{eq4.12}) - \eqref{eq4.13} has a unique optimal solution $(u^*,f^*) \in \mathcal{U}_{ad}$ for almost surely $\omega \in \Omega$.
\end{lemma}

\subsubsection{Discretization of the optimality system} \label{subsubsec4.2.1}

In this article, we adopt the Discretize-then-Optimize strategy, which is a popular and powerful method to solve the optimal control problem \citep{chen2022reduced, leykekhman2012investigation, liu2019non, mathew2007analysis}. Similarly to Section \ref{subsubsec4.1.2}, we make the spatial and discretization of the distributed optimal control problem in \eqref{eq4.12}) - \eqref{eq4.13} using the standard FEM method. Let $\{\phi_j\}_{j=1}^N$ denote the finite element basis functions. The state and control variables are expressed as below: 

$$u_h (\bm{x},t,\omega) = \sum_{j=1}^N u_j(t,\omega) \phi_j(\bm{x})\quad \text{and} \quad f_h (\bm{x},t) = \sum_{j=1}^N f_j(t) \phi_j(\bm{x}).$$

Given $M$ i.i.d. samples $\{\omega_m\}_{m=1}^M$ randomly drawn from $\Omega$ by the probability measure $\mathbf{P}$, the state equations in \eqref{eq4.13} are discretized as follows:

\begin{equation} 
\begin{aligned}
\sum_{j=1}^N \frac{\mathrm{d}u_{j,m}}{\mathrm{d}t} \int_D \phi_j \phi_i d\bm{x} \:&+\:\sum_{j=1}^N u_{j,m} \left[\int_D \bar{a} \nabla \phi_j \nabla \phi_i d\bm{x}  + \:  \int_D \widetilde{a}_m \nabla \phi_j \nabla \phi_i d\bm{x} \right]\\
& \: = \: \sum_{j=1}^N f_j \int_D \phi_j \phi_i d\bm{x}, \quad i = 1,\cdots,N, m = 1,\cdots,M.
\nonumber
\end{aligned}
\end{equation}

\noindent which can be also expressed as the following algebraic system:

\begin{equation} 
\mathbf{G} \frac{\mathrm{d}}{\mathrm{d}t}\bm{u}_m(t) + \left(\bar{\mathbf{A}} \:+\: \widetilde{\mathbf{A}}_m\right)\bm{u}_m(t) \:=\: \mathbf{G}\, \bm{f}(t), \quad t \in (0,T]. 
\nonumber
\end{equation}

\noindent The mass matrix $\mathbf{G}$, the stiffness matrices $\bar{\mathbf{A}}$ and $\{\widetilde{\mathbf{A}}_m\}_{m=1}^M$ are defined in \eqref{eq4.6}. The solutions are given by

$$\bm{u}_m(t) = \left( u_{1,m}(t), ..., u_{N,m}(t)\right)^T \quad \text{and}\quad \bm{f}(t) = \left( f_1(t), ..., f_N(t)\right) ^T.$$

\noindent The semi-discrete objective functional is stated as below:
	
$$\hat{\mathcal{J}}(\bm{u}_m, \bm{f}) \: := \: \frac{1}{M} \sum_{m=1}^M \left[\frac{1}{2} \int_0^T \left(\bm{u}_m(t) - \bm{U}(t)\right)^T \mathbf{G} \left(\bm{u}_m(t) - \bm{U}(t)\right) \mathrm{d}t\right] + \frac{\beta}{2} \int_0^T \bm{f}(t)^T\mathbf{G}\bm{f}(t) \mathrm{d}t ,$$

\noindent where the desired state is $(\bm{U}(t))_i = \int_\mathcal{D} U(\bm{x}, t) \phi_i(\bm{x}) \mathrm{d}\bm{x}$.

Here we adopt the backward Euler scheme to deal with the temporal space. Denote $\bm{u}_{m,l} := \bm{u}_m(t_l)$, $\bm{f}_{l} := \bm{f}(t_l)$ and $\bm{U}_{l} := \bm{U}(t_l)$ at each time step $l = 0,1,\cdots,L$. The the fully discrete objective functional is given by

\begin{equation}\label{eq4.15}
\hat{\mathcal{J}}(\bm{u}_{m,l}, \bm{f}_l) \: := \: \frac{1}{M} \sum_{m=1}^M \sum_{l=1}^L \frac{\Delta t}{2} \left(\bm{u}_{m,l} - \bm{U}_l\right)^T \mathbf{G} \left(\bm{u}_{m,l} - \bm{U}_l\right)  + \sum_{l=1}^L \frac{\beta \Delta t}{2} \bm{f}_l^T\mathbf{G}\bm{f}_l,
\end{equation}

\noindent where $\bm{u}_{m,l}$ and $\bm{f}_l$ are the numerical solutions of the fully discrete state equations as below:

\begin{equation}
\mathbf{G}\frac{\bm{u}_{m,l+1}-\bm{u}_{m,l}}{\Delta t} + \left(\bar{\mathbf{A}} \:+\: \widetilde{\mathbf{A}}_m\right)\bm{u}_{m,l+1} \:=\: \mathbf{G}\, \bm{f}_{l+1}.
\nonumber
\end{equation}

Finally, we can derive the following algebraic system arising from the discrete approximation of the constraint unsteady diffusion SPDE:

\begin{equation}\label{eq4.16}
\left\{\:
\begin{aligned}
& \left(\mathbf{G} + \Delta t (\bar{\mathbf{A}} + \widetilde{\mathbf{A}}_m)\right)\bm{u}_{m,l+1} \:=\: \mathbf{G}\,\bm{u}_{m,l}\:+\:\Delta t \,\mathbf{G} \,\bm{f}_{l+1},\\
& \,\bm{u}(0;\omega_m) \:=\: \bm{u}_0,
\end{aligned}
\right.
\end{equation}

\noindent for $m=1,\cdots,M$ and $l = 0,1,\cdots,L-1$. The fully discretized optimal control problem is then given as follows: minimize the cost functional $\hat{\mathcal{J}}(\bm{u}_{m,l}, \bm{f}_l)$ in \eqref{eq4.15}, subject to the $M$ constrained equations in \eqref{eq4.16}. Then, our next step is to apply the LRNS solver to solve the resulting linear systems efficiently.

\subsubsection{Construct the LRNS solver for the unsteady stochastic diffusion optimal control problem}\label{subsubsec4.2.2}

To efficiently solve the optimality system in \eqref{eq4.15}-\eqref{eq4.16}, we introduce the so-called reduced approach \cite{borzi2009multigrid, chen2022reduced, meidner2008priori, zhu2024low}, which aims to eliminate the SPDE constraints and construct an unconstrained optimization problem. First, we rewrite the discrete state equations in \eqref{eq4.16}, and the approximate state for $l$-th time step of the $m$th Monte-Carlo realization is given by

\begin{equation}\label{eq4.17}
\bm{u}_{m,l+1} \enspace=\enspace \Delta t \left(\mathbf{G} + \Delta t (\bar{\mathbf{A}} + \widetilde{\mathbf{A}}_m)\right)^{-1}\,\mathbf{G} \,\bm{f}_{l+1}\:+\:\left(\mathbf{G} + \Delta t (\bar{\mathbf{A}} +\widetilde{\mathbf{A}}_m)\right)^{-1}\mathbf{G}\,\bm{u}_{m,l}.
\end{equation}

To reduce the space and computational complexities from inverting these large-scale coefficient matrices, we construct the LRNS solver for the distributed optimal control problem. We begin by introducing the following matrices:

\begin{equation}\label{eq4.18}
\bar{\mathbf{K}} \::=\: \frac{\mathbf{G}}{\Delta t} + \bar{\mathbf{A}}, \quad \widetilde{\mathbf{K}}_m \::=\: \widetilde{\mathbf{A}}_m,\qquad m = 1,\cdots,M.
\end{equation}

The perturbation matrices are approximated by the low-rank representation $\widetilde{\mathbf{K}}_m \approx \mathbf{U} \mathbf{V}^T$ obtained by Algorithm \ref{alg1}. Then, we deal with the inverse matrices in \eqref{eq4.17} by the Neumann matrix series as below:

\begin{equation}
\Delta t \left(\mathbf{G} + \Delta t (\bar{\mathbf{A}} + \widetilde{\mathbf{A}}_m)\right)^{-1} \:=\: \left(\bar{\mathbf{K}} + \widetilde{\mathbf{K}}_m\right)^{-1} \:\approx\: \sum_{r=0}^R \left(-\bar{\mathbf{K}}^{-1}\mathbf{U}\mathbf{V}_m^T\right)^{r} \bar{\mathbf{K}}^{-1}.
\nonumber
\end{equation}

\noindent Denote $\mathbf{Z}_m := \sum_{r=0}^R \left(-\bar{\mathbf{K}}^{-1}\mathbf{U}\mathbf{V}_m^T\right)^{r} \bar{\mathbf{K}}^{-1} \mathbf{G}$ for $m = 1,\cdots,M$. Then, the linear equation in \eqref{eq4.17} reduces to:

\begin{equation}
\bm{u}_{m,l+1} \enspace=\enspace \mathbf{Z}_m\,\bm{f}_{l+1}\:+\:\frac{1}{\Delta t}\,\mathbf{Z}_m\bm{u}_{m,l}, \qquad m = 1,\cdots,M,\:l = 0,\cdots,L-1.
\nonumber
\end{equation}

For notational simplicity, let $\bm{u}_{m,l+1} := \mathcal{U}_{m,l}(\bm{f}_{l+1})$ denote the control-to-state mapping, which presents the dependence of $\bm{u}_{m,l+1}$ on $\bm{f}_{l+1}$. In this case, there exists an unique state $\bm{u}_{m,l+1}$ for a given control $\bm{f}_{l+1}$. Then, we plug the mapping into the objective functional in \eqref{eq4.15} and introduce the following reduced cost functional:

\begin{equation}\label{eq4.19}
\hat{\mathcal{J}}( \bm{f}_l) \: = \: \frac{1}{M} \sum_{m=1}^M \sum_{l=1}^L \frac{\Delta t}{2} \left(\mathcal{U}_{m,l-1}(\bm{f}_{l}) - \bm{U}_l\right)^T \mathbf{G} \left(\mathcal{U}_{m,l-1}(\bm{f}_{l}) - \bm{U}_l\right)  + \sum_{l=1}^L \frac{\beta \Delta t}{2} \bm{f}_l^T\mathbf{G}\bm{f}_l.
\end{equation}

Therefore, we transform the original stochastic optimal control problem in \eqref{eq4.1}-\eqref{eq4.2} into the unconstrained optimization problem in \eqref{eq4.19}, which can significantly reduce computational cost and allow the use of various gradient-based algorithms, such as the steepest descent method and the Newton’s method. Based on the chain rule, it is straightforward to derive the gradient of the cost functional:

\begin{equation} \label{eq4.20}
\nabla \hat{\mathcal{J}}(\bm{f}_l) \: = \: \frac{\mathrm{d} \hat{\mathcal{J}}}{\mathrm{d} \bm{f}_l} \: = \: \frac{1}{M} \sum_{m=1}^M \left( \frac{\partial \hat{\mathcal{J}}}{\partial \bm{u}_{m,l}} \right)^T \frac{\mathrm{d} \bm{u}_{m,l}}{\mathrm{d} \bm{f}_l}  + \frac{\partial \hat{\mathcal{J}}}{\partial \bm{f}_l},
\end{equation}

\noindent where the partial derivatives $\frac{\partial \hat{\mathcal{J}}}{\partial \bm{u}_{m,l}}$, $\frac{\partial \hat{\mathcal{J}}}{\partial \bm{f}_l}$, and the sensitivity $\frac{\mathrm{d} \bm{u}_{m,l}}{\mathrm{d} \bm{f}_l}$ are respectively given by

\begin{equation} \label{eq4.21}
\frac{\partial \hat{\mathcal{J}}}{\partial \bm{u}_{m,l}} \: = \: \Delta t\,\mathbf{G} \left( \mathcal{U}_{m,l-1}(\bm{f}_l) - \bm{U}_l\right) , \qquad \frac{\partial \widehat{\mathcal{J}}}{\partial \bm{f}_h} \: = \: \beta \Delta t\,\mathbf{G} \,\bm{f}_l, \qquad \frac{\mathrm{d} \bm{u}_{m,l}}{\mathrm{d} \bm{f}_l} \: = \: \mathbf{Z}_m.
\end{equation}

By substituting \eqref{eq4.21} into \eqref{eq4.20}, we obtain the following gradient of the objective functional:

\begin{equation} \label{eq4.22}
\nabla \hat{\mathcal{J}}(\bm{f}_l) \: = \: \frac{1}{M} \sum_{m=1}^M \Delta t \,\mathbf{Z}_m^T \mathbf{G} \left( \mathbf{Z}_m \bm{f}_l - \bm{U}_l\right) \,+\, \beta \Delta t \,\mathbf{G} \bm{f}_l.
\end{equation}

\noindent Similarly, we can derive the following Hessian matrix of of the objective functional:

\begin{equation} \label{eq4.23}
\Delta \hat{\mathcal{J}} \: = \: \frac{1}{M} \sum_{m=1}^M \Delta t \,\mathbf{Z}_m^T \mathbf{G}  \mathbf{Z}_m \,+\, \beta \Delta t \,\mathbf{G}.
\end{equation}

Finally, the control variable is updated in the unconstrained optimization procedure by

\begin{equation} \label{eq4.24}
\bm{f}_l^{(k+1)}\: = \: \bm{f}_l^{(k)}  \; + \; \delta \bm{f}_l^{(k)} \: = \: \bm{f}_l^{(k)}  \; - \; \alpha^{(k)} \mathbf{d}_l^{(k)},
\end{equation}

\noindent where $k$ denotes the iteration index and $\alpha^{(k)} > 0$ is the step length. The descent direction $\mathbf{d}_l^{(k)}$ can be determined based on the first-order and second-order derivatives in \eqref{eq4.20}-\eqref{eq4.21}. It is noteworthy that the Hessian matrix $\Delta \hat{\mathcal{J}}$ is independent of the control variable $\bm{f}_l$, i.e., it remains constant throughout the iterations and does not require recomputation, which significantly reduces the computational cost. The pseudo-code of the LRNS solver for the distributed optimal control problem governed by the unsteady stochastic diffusion equation is summarized in Algorithm \ref{alg4}.

\begin{algorithm}
\caption{The LRNS solver for the distributed optimal control problem in \eqref{eq4.12} governed by the unsteady stochastic diffusion equation with random permeability in \eqref{eq4.13}} \label{alg4}
\begin{algorithmic}[1]
\Require
A tessellation $\mathcal{T}_h$ of spatial domain $\mathcal{D}$, final time $T$, MC sample size $M$, numbers of spatial and temporal discretization for SPDEs $N$ and $L$, Neumann series truncation index $R$, and data compression ratio $\tau$, optimization tolerance $\varepsilon$.
\Ensure
Numerical optimal control $\{\bm{f}_l^{*}\}_{l=0}^L$ and , and the expectation of optimal state $\{\mu_l^{*}\}_{l=0}^L$.
			
\State Construct the finite element space $V_h \in \mathcal{H}^1(\mathcal{D})$.

\State Generate the deterministic and perturbation permeability $\bar{a},\{\widetilde{a}_m\}_{m = 1}^M$ by the truncated KL expansion in \eqref{eq4.3}.
		
\State Assemble the stiffness matrices $\bar{\mathbf{K}}, \{\mathbf{\widetilde{K}}_m\}_{m=1}^M$ by \eqref{eq4.18}, and deal with the Dirichlet boundary conditions of $\bar{\mathbf{K}}$.
			
\State Determine the reduced dimension $k = \lceil \tau N \rceil$ and compute the optimal rank-$k$ approximation $\mathbf{U}, \{\mathbf{V}_m\}_{m=1}^M$ of $\{\widetilde{\mathbf{K}}_m\}_{m=1}^M$ using Algorithm \ref{alg1}.

\State Compute the matrice $\mathbf{Z}_m := \sum_{r=0}^R \left(-\bar{\mathbf{K}}^{-1}\mathbf{U}\mathbf{V}_m^T\right)^{r} \bar{\mathbf{K}}^{-1} \mathbf{G}$.

\State Generate initial guess of control $\bm{f}_l^{(0)}$ and solve $\bm{u}_{m,l}^{(0)}=\mathcal{U}_{m,l-1}(\bm{f}_l^{(0)})$.
			
\While{ $\Vert \nabla \widehat{\mathcal{J}}(\bm{f}_l^{(k)}) \Vert \: > \: \varepsilon$ }
			
\State Compute the gradient $\nabla \widehat{\mathcal{J}}(\bm{f}_l^{(k)})$ or the hessian matrix $\Delta \widehat{\mathcal{J}}$ with respect to $\bm{f}_l^{(k)}$ by \eqref{eq4.23} and \eqref{eq4.24}.
			
\State Choose a suitable step size $\alpha^{(k)}  >  0$.
			
\State Compute the step $\delta \bm{f}_l^{(k)}$ based on the derivative information $\nabla \widehat{\mathcal{J}}(\bm{f}_l^{(k)})$ or $\Delta \widehat{\mathcal{J}}$, and the step size $\alpha^{(k)}$.
			
\State Update the control variable $\bm{f}_l^{(k+1)}$ by \eqref{eq4.24}.
			
\State $k = k + 1$.
			
\EndWhile
			
\State \Return $\bm{f}_l^* = \bm{f}_l^{(k)}$ and $\mu_l^* = \frac{1}{M} \sum_{m=1}^M \bm{u}_{m,l}^{(k)}$ for $m = 1,\cdots,M$ and $l = 0,\cdots,L$.
			
\end{algorithmic}
\end{algorithm}

\subsubsection{Numerical experiments on the unsteady stochastic diffusion optimal control problem} \label{subsubsec4.2.3}

\begin{figure}[ht]
\centering
\includegraphics[width=1 \textwidth]{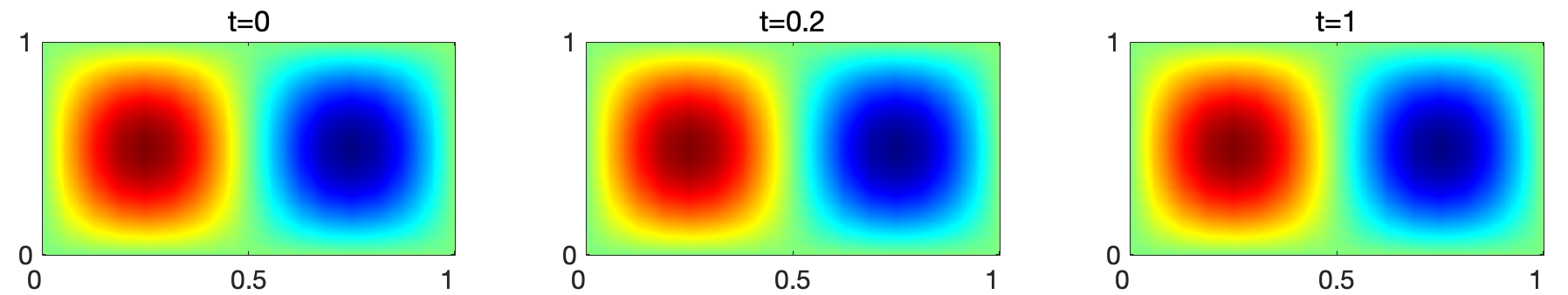}
\caption{Desired state at $t=0,0.2,1$.}
\label{fig4.8}
\end{figure}

In the following, we demonstrate the numerical performance of the LRNS solver to the distributed optimal control problem in \eqref{eq4.12} governed by the unsteady stochastic diffusion equation SPDE in \eqref{eq4.13}. The numerical settings are similar to Section \ref{subsubsec4.1.4}. We define the regularization parameter $\beta = 1 \times 10^{-3}$, the tolerance parameter $\epsilon = 1 \times 10^{-3}$, and the maximum iterations count in the line search method with the Wolfe condition $it_{max} = 50$ for the numerical implementation. The initial condition is set to be $u_0(x,y) = \sin(2\pi x)\sin(\pi y)$. The distributed optimal control problem has the desired state $U(\bm{x})$ $U(x,y,t) = e^{-\pi t}\sin(2 \pi x) \sin(\pi y)$, which is visualized in Figure \ref{fig4.8}.

\begin{figure}[ht]
\centering
\includegraphics[width=1 \textwidth]{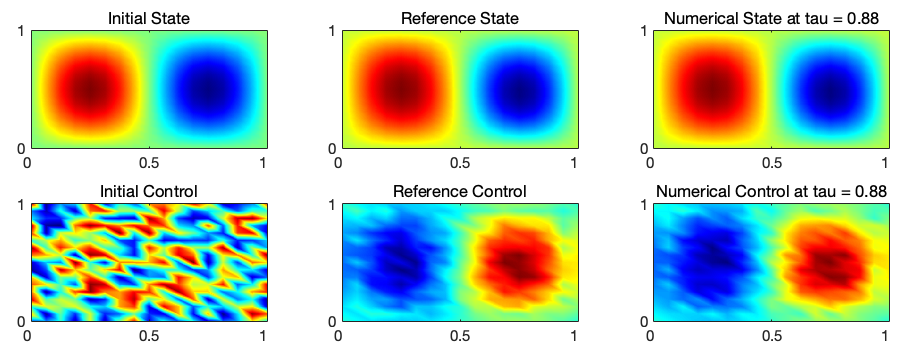}
\caption{The initial state $u^{(0)}(\bm{x})$ and control $f^{(0)}(\bm{x})$ (left), the reference state $u_{ref}(\bm{x})$ and control $f_{ref}(\bm{x})$ (middle), and the numerical optimal state $u^*(\bm{x})$ and control $f^*(\bm{x})$ using the data compression ratio $\tau =88\%$ (right).}
\label{fig4.9}
\end{figure}

We first focus on the efficacy of the proposed solver under a representative data compression level $\tau = 88\%$. The reference solution is still computed using the classical MC-FEM method, and the resulting optimization problem is solved by Newton’s method. As shown in Figure \ref{fig4.9}, noticeable changes are observed in both the state and control variables after optimization. The numerical solutions obtained by the LRNS solver closely match the reference solution, which confirms the feasibility and effectiveness of Algorithm \ref{alg4}.
	
Then we investigate the sensitivity of the proposed algorithm to the data compression ratio $\tau$. Table \ref{tab4.2} compares the optimization results achieved by the LRNS solver under different choices of $\tau = 100\%, 95\%, 88\%, 50\%, 30\%, 10\%$. A relatively large $\tau$ leads to better optimization efficacy and we can achieve a small ratio between the initial and final objective values $\widehat{J}^{*}/\widehat{J}^{0}$. In contrast, selecting a small $\tau$ significantly reduces the CPU run-time. These observations are consistent witjh our conclusions in Table \ref{tab3.1} and Theorem \ref{th3.7}.

\begin{table}[ht]
\centering 
\begin{tabular}{cccccccc}
\hline
\noalign{\vskip 2pt}
Data compression ratio & Tic-toc time (s) & $\nabla \widehat{J}^{*}$ & $\widehat{J}^{(0)}$ & $\widehat{J}^{*}$ & $\widehat{J}^{*}/\widehat{J}^{0}$  \\ 
\noalign{\vskip 2pt}\hline
$100\%$ & $5843.9298$ & $9.3012 \times 10^{-4}$ & $0.06579257$ & $0.00142084$ & $2.16\%$\\
$95\%$  & $5566.0615$ & $8.9203 \times 10^{-4}$ & $0.06579152$ & $0.00143315$ & $2.18\%$\\
$88\%$  & $4655.2956$ & $8.7880 \times 10^{-4}$ & $0.06579140$ & $0.00145132$ & $2.21\%$\\
$50\%$  & $4362.0363$ & $8.7340 \times 10^{-4}$ & $0.06702326$ & $0.00214603$ & $3.20\%$\\
$30\%$  & $3613.6498$ & $8.3761 \times 10^{-4}$ & $0.06712217$ & $0.00217054$ & $3.23\%$\\
$10\%$  & $3169.0377$ & $8.3048 \times 10^{-4}$ & $0.06732181$ & $0.00225543$ & $3.35\%$\\\hline
\end{tabular}
\caption{Simulation results for Algorithm \ref{alg4} about the tic-toc time, final gradient value $\nabla\widehat{J}^{*}$, initial and final objective values $\widehat{J}^{(0)}$ and $\widehat{J}^{*}$, and the ratio $\widehat{J}^{*}/\widehat{J}^{(0)}$ under six different data compression ratios $\tau = 100\%, 95\%, 88\%, 50\%, 30\%, 10\%$.}
\label{tab4.2}
\end{table}

According to the reduced formulation \eqref{eq4.19}, the original distributed optimal control system is reformulated as an unconstrained optimization problem. Therefore, an appealing property of Algorithm \ref{alg4} is its flexibility to incorporate different kinds of gradient-based algorithms. To illustrate this advantage, we employ three widely used unconstrained optimization algorithms within the LRNS framework: the steepest descent method, the stochastic gradient descent method and the Newton's method. The results in Table \ref{tab4.3} show that all three algorithms achieve satisfactory optimization performance, which validates the effectiveness and robustness of Algorithm \ref{alg4}.

\begin{table}[ht]
\centering 
\begin{tabular}{cccccc}
\hline
\noalign{\vskip 2pt}
Data compression ratio & Iteration & Tic-toc time (s) & $\widehat{J}^{(0)}$ & $\widehat{J}^{*}$ & $\widehat{J}^{*}/\widehat{J}^{0}$  \\ 
\noalign{\vskip 2pt}\hline
Newton's Method & $48$ & $4655.2956$ & $0.06579140$ & $0.00145132$ & $2.21\%$\\
Steepest Descent Method & $96$ & $8450.5518$ & $0.06579140$ & $0.00149239$ & $2.27\%$\\
Stochastic Gradient Method & $376$ & $3815.8394$ & $0.06579140$ & $0.00151902$ & $2.31\%$\\\hline
\end{tabular}
\caption{Simulation results for three different unconstrained optimization methods about the iteration, tic-toc time, initial and final objective values $\widehat{J}^{(0)}$ and $\widehat{J}^{*}$, and the ratio $\widehat{J}^{*}/\widehat{J}^{(0)}$.}
\label{tab4.3}
\end{table}

\section{Conclusions}\label{sec5}

In this article, we develop an efficient numerical solver for unsteady diffusion-type partial differential equations with random coefficients. The proposed solver combines a novel generalized low-rank matrix approximation with the Neumann matrix series expansion to reduce space and computational complexities without losing numerical accuracy. We apply it to address two representative applications: unsteady stochastic diffusion equations and distributed optimal control problems governed by unsteady diffusion SPDE constraints. The error analysis is carried out and numerical experiments are conducted to validate the proposed solver and the theoretical conclusions.

Several directions for future work can be pursued. One natural extension is to apply the proposed framework to nonlinear and unsteady stochastic partial differential equations, such as the time-dependent Navier–Stokes equations. Another promising direction is to integrate the solver with advanced sampling strategies and adaptive meshing techniques. These extensions would further broaden the applicability of the proposed approach in UQ.

\section*{Acknowledgments}
The authors would like to thank the anonymous referees and the editor for their valuable comments and suggestions, which led to considerable improvement of the article.

\section*{Conflict of Interest}
All authors declare that there are no conflicts of interest regarding the publication of this paper.

\bibliographystyle{elsarticle-harv}
\bibliography{bib}

\end{document}